\documentclass[11pt]{article}

\usepackage{amsthm}

\usepackage{epsfig}

\usepackage{latexsym}
\usepackage{amsmath}
\usepackage{amssymb}
\usepackage{amsfonts}

\usepackage[size=1.5em,nohug,heads=LaTeX,PostScript=dvips]{diagrams}

\theoremstyle{definition}

\newtheorem{dfn}{Definition}[section]
\newtheorem{defn}[dfn]{Definition}
\newtheorem{ex}[dfn]{Example}
\newtheorem{rem}[dfn]{Remark}

\theoremstyle{plain}

\newtheorem{thm}[dfn]{Theorem}
\newtheorem{example}[dfn]{Example}

\newtheorem{lem}[dfn]{Lemma}

\newtheorem{prop}[dfn]{Proposition}
\newtheorem{prob}[dfn]{Problem}

\newtheorem{cor}[dfn]{Corollary}

\newtheorem{convention}[dfn]{Convention}
\newtheorem{notation}[dfn]{Notation}

\def\proof{\par\medskip\noindent{\it Proof. }}

\def\C{{\mathbb C}}
\def\R{{\mathbb R}}

\def\Z{{\mathbb Z}}

\def\d{{\mathcal D}}

\def\N{{\mathbb N}}

\def\c{{\mathcal C}}
\def\e{{\mathcal E}}

\def\d{{\mathcal D}}

\def\eps{\epsilon}
\def\al{\alpha}
\def\be{\beta}
\def\ga{\gamma}
\def\Ga{\Gamma}

\def\del{\delta}
\def\De{\Delta}
\def\Del{\Delta}
\def\Si{\Sigma}
\def\si{\sigma}

\def\la{\lambda}
\def\La{\Lambda}
\def\Om{\Omega}
\def\om{\omega}

\def\acts{\curvearrowright}
\def\D{\partial}

\def\embed{\hookrightarrow}
\def\8{\infty}
\def\<{\langle}
\def\>{\rangle}
\def\geo{\partial_{\infty}}

\def\ol{\overline}

\def\t{\tilde}

\newcommand{\diam}{\operatorname{diam}}

\begin{document}

\title{Affine buildings for dihedral groups}
\author{Arkady Berenstein and Michael Kapovich}
\date{September 1, 2008}

\maketitle

\begin{abstract}
We construct rank 2 thick nondiscrete affine buildings associated
with an arbitrary finite dihedral group.
\end{abstract}

\section{Introduction}

In his foundational work on buildings, J.\,Tits \cite{Tits1,
Tits3} proved that all thick irreducible spherical and affine
buildings of rank $\ge 3$ (including nondiscrete ones) are
associated with algebraic groups and, hence, have Weyl groups
coming from root systems. He then showed \cite{Tits2} how to use
the {\em free construction} to prove existence of thick spherical
buildings of rank 2 modeled on arbitrary finite Coxeter groups
(see Definition \ref{modeled}). Such buildings are necessarily
infinite except for dihedral groups of the order $2, 4, 6, 8, 12,
16$. M.\,Ronan \cite{Ronan1} constructed a variety of irreducible
rank 2 discrete affine buildings not coming from algebraic groups.
The corresponding finite Coxeter groups are necessarily of
crystallographic type, i.e., they are the dihedral groups $D_m,
m=3, 4, 6$. However, examples of nondiscrete affine buildings of
rank 2 corresponding to the rest of finite dihedral groups
remained elusive. The main result of this paper is a construction
of such buildings.

Consider the dihedral group $W=D_m$ of order $2m, m\ge 2$. Take
 a countable dense subgroup $\La$ of translations of the apartment
 $A=\R^2$, invariant under the action of $W$ by conjugations, so that for
 the group $W_{af}=\La \rtimes W$ (which is, in fact, a Coxeter group),
 every vertex of the Coxeter complex $(A,W_{af})$ is special. Such Coxeter complexes will be
called {\em special}. For instance, when $m$ is a prime number,
the complex $(A,W_{af})$ is special no matter what $\La$ is. The
main goal of this paper is to prove

\begin{thm}\label{main}
For every special Coxeter complex $(A,W_{af})$ with nondiscrete
countable group $W_{af}$, there exists a thick Euclidean building
$X$ of rank 2 modeled on $(A,W_{af})$.
\end{thm}

After proving this theorem we found out that for $m=8$, existence
of thick Euclidean buildings (for some countable $W_{af}$) was
also  proven in a recent work of P.Hitzelberger, L.Kramer and
R.Weiss \cite{HKW} by a totally different method: Their buildings
are embedded in  buildings of the type $F_4$.

It is quite likely that Theorem \ref{main} holds for semidirect
products $W_{af}$ of $W$ with arbitrary (nondiscrete) groups of
translations $\La\subset \R^2$, or, at least countable ones.

\medskip
The strategy for proving Theorem \ref{main} is the same as in
Tits' free construction of rank 2 spherical buildings, see
\cite{Tits2, FS}, i.e., we run a {\em ``Ponzi-scheme''}:

We inductively construct a sequence of CAT(0) spaces
$$
X_1\to X_2 \to ... \to X_n \to X_{n+1} \to ...
$$
where each $X_n$ is {\em finite}, i.e., is a finite union of
apartments. (The embeddings $X_n\to X_{n+1}$ are {\em not}
isometric.) Each $X_n$ is modeled on $(A, W_{af})$ but fails,
rather badly, the rest of axioms of a thick building, as ``most''
pairs of points are not contained in a common apartment. Given
points (say, vertices) $x_1, x_2\in X_n$, we obtain $X_{n+1}$ by
attaching an apartment to $X_n$ along a closed convex subcomplex
$C\subset X_n$ containing $x_1, x_2$.  Then $x_1, x_2\in X_{n+1}$
belong to a common apartment. This, of course, produces more
points in $X_{n+1}$ which do not belong to a common apartment, as
well as more walls which are not {\em thick}. However, these will
be taken care of at some later stages of the construction. The
hardest part of the proof is to establish the following {\em
combinatorial convexity theorem}, proven in sections
\ref{convexity} and \ref{generalized}:

\begin{thm}\label{conv}
Let $X$ be a CAT(0) metric space, which is a finite space modeled
on $(A, W_{af})$. Then, given any pair of points $x_1, x_2\in X$
and germs of Weyl chambers $\si_i\subset X$ with the tips at
$x_i$, there exists a {\em planar} convex subcomplex $C\subset X$
containing $x_i$ and $\si_i, i=1, 2$. The same applies to each
pair of {\em super-antipodal} chambers $\si_1, \si_2$ at infinity,
i.e., pair of chambers which admit  regular points connected by a
geodesic in $X$.
\end{thm}

Planarity of $C$ means existence of a {\em weakly isometric
embedding} $C\embed A$. (We introduce $D$-isometry and weak
isometry in Sections \ref{metric} and \ref{spaces} as certain
generalizations of the ordinary isometry between metric spaces.)
The embedding $C\embed A$ allows one to attach  $A$ to $X$ along
$C$ so that $X\cup_C A$ is again a CAT(0) space which is finite,
modeled on $(A, W_{af})$. The convex subcomplexes $C$ used in our
paper are (extended) {\em corridors} obtained by taking a finite
union of {\em tiles} in $X$. The appearance of corridors in this
context is natural since, in the case when $X$ is a Euclidean
building, every pair of vertices $x_1, x_2\in X$ admits a
combinatorial convex hull, called a {\em Weyl hull}, which is a
tile with the vertices $x_1, x_2$.

Since each $X_n$ will have only countably many corridors, we
enumerate them and form the sequence $(X_n)$ so that for every
convex corridor $C$ in $X_k$ there exists $n\ge k$ such that
$X_{n+1}$ is formed by attaching an apartment along a convex
corridor $C'\subset X_n$ containing $C$.

Organization of the paper. In sections \ref{complexes},
\ref{metric}, \ref{spaces} we review standard concepts related to
buildings: Coxeter complexes, CAT(0) and CAT(1) metric spaces and
spaces modeled on Coxeter complexes. We review the building axioms
in section \ref{spaces}. The only nonstandard material in these
sections is the one of {\em $D$-isometries} between metric spaces
and the corresponding concept of {\em weak isometries} in the
context of spaces modeled on $(A, W_{af})$. In section
\ref{corridors} we define the main technical tool of this paper
--- {\em corridors}, which are finite unions of {\em tiles}, parallelograms in $A$ with
the angles $\frac{\pi}{m}$ and $\pi-\frac{\pi}{m}$. In section
\ref{orders} we define two orders on corridors and prove existence
of maximal elements with respect to these orders. The {\em maximal
corridors} are proven to be convex in section \ref{convexity}. In
the same section we prove a weaker version of Theorem \ref{conv},
i.e., we prove that every two points in $X$ belong to a corridor
(and, therefore, to a convex corridor). In section
\ref{generalized} we introduce the notion of {\em extended
corridors}, which are noncompact analogues of corridors and finish
the proof of Theorem \ref{conv}. In section \ref{category} we
define {\em enlargement} of a space $X$ modeled on $(A, W_{af})$,
which is the operation of attaching to $X$ an apartment along an
extended convex corridor in $X$. We then show that the image of a
convex corridor under enlargement is again a corridor (possibly
nonconvex). In section \ref{mainproof} we put it all together and
prove our main result, Theorem \ref{main}, by constructing an
increasing sequence $(X_n)$ of spaces modeled on $(A, W_{af})$,
whose union is an affine building. In section \ref{conclusion} we
make some speculative remarks on possible generalizations of our
results, in particular, constructing rank 2 buildings with highly
transitive automorphism groups.

\medskip
{\bf Acknowledgments.} Our collaboration on this project started
at the AIM workshop ``Buildings and Combinatorial Representation
Theory'' in 2007 and we are grateful to AIM for this opportunity.
The first author was supported by the NSF grant DMS-08-00247. The
second author was supported by the NSF grant DMS-05-54349. He is
grateful to Linus Kramer and Petra Hitzelberger for the reference
to their work.

\tableofcontents

\section{Coxeter complexes}\label{complexes}

Let $A$, the {\em apartment}, be either the Euclidean space $E^n$
or the unit $n-1$-sphere $S^{n-1}\subset E^n$ (we will be
primarily interested in the case of Euclidean plane and the
circle). If $A=S^{n-1}$, a {\em Coxeter group} acting on $A$ is a
finite group $W$ generated by isometric reflections. If $A=E^n$, a
{\em Coxeter group} acting on $A$ is a group $W_{af}$ generated by
isometric reflections in hyperplanes in $A$, so that the linear
part of $W_{af}$ is a Coxeter group acting on $S^{n-1}$. Thus,
$W_{af}=\La \rtimes W$, where $\La$ is a certain (countable or
uncountable) group of translations in $E^n$.

\begin{defn}
A spherical or Euclidean\footnote{Also called {\em affine}.}
{\em Coxeter complex} is a pair $(A, G)$, of the form $(S^{n-1}, W)$ or $(E^n, W_{af})$.
\end{defn}

A  {\em wall} in the Coxeter complex $(A, G)$ is the fixed-point
set of a reflection in $G$. A {\em regular point} in a Coxeter
complex is a point which does not belong to any wall. A {\em
singular point} is a point which is not regular. A {\em
half-apartment} in $A$ is a closed half-space bounded by a wall. A
{\em cell} in is a convex polyhedron which is the intersection of
finitely-many half-apartments.

\begin{rem}
Note that in  the spherical case, there is a natural cell complex
in $S^n$ associated with $W$, where the {\em cells} are those
cells (in the above sense) which cannot be separated by a wall.
However, the affine case, when $W_{af}$ is nondiscrete, there will
be no cell natural complex attached to $W_{af}$.
\end{rem}

{\em Chambers} in $(S^{n-1}, W)$ are the fundamental domains for
the action $W\acts S^{n-1}$, i.e., the closures of the connected
components of the complement to the union of walls.

A {\em dilation} of an affine Coxeter complex is a map of the form $x\mapsto \la x+ v$, where
$v$ is a vector and $\la>0$. Thus, we regard translations as a (limiting form of) dilations.

A {\em vertex} in $(A, G)$ is a (component of, in the spherical
case) the 0-dimensional intersection of walls. We will consider
almost exclusively only those Coxeter complexes which have at
least one vertex; such complexes are called {\em essential}.
Equivalently, these are spherical complexes where the group $G$
does not have a global fixed point and those Euclidean Coxeter
complexes where $W$ does not have a fixed point in $S^{n-1}$.

A {\em special point} in $(E^n, W_{af})$ is a point whose
stabilizer in $W_{af}$ is isomorphic to $W$. We will always choose
coordinates in $A=E^n$ so that the origin is special. A Coxeter
complex where every vertex is special is called a {\em special
Coxeter complex}. In general, the stabilizer of a point $x\in A$
is naturally isomorphic a Coxeter subgroup of $W$.

Consider now the case of $n=2$, when $W=D_m$, the dihedral group
of order $2m$. Under the identification of the vertices of the
regular $2m$-gon with the $2m$-th roots of unity, the group $\La$
always contains (or equals to) the ring $\La_0=\Z[\zeta]$, where
$\zeta$ is a primitive $2m$-th root of unity if $m$ is odd and a
primitive $m$-th root of unity if $m$ is even. It is well--known
that $\La_0$ is discrete in $\C=\R^2$ if and only if $m=1, 2, 3,
4, 6$.

If $m$ is prime, then every proper Coxeter subgroup of $W$ is
isomorphic to $\Z/2$ and, hence, every vertex of $(E^2, W_{af})$
is special. If $m$ is not prime and $\La$ is the smallest
nontrivial group of translations normalized by $W$, the complex
$(A, W_{af})$ is non-special. Nevertheless,

\begin{lem}
For every (countable) $W_{af}=\La\rtimes W$, there exists a
(countable) group $\La'$ containing $\La$ and normalized by $W$, so
that the Coxeter complex $(A, \La'\rtimes W)$ is special.
\end{lem}
\proof Set $\La^1:=\La$ and $W^1_{af}:= W_{af}$ and $C^1:= (A,
W^1_{af})$. We continue inductively: Given a Coxeter complex
$C^i:= (A, W^i_{af})$ with (countable) Coxeter group, we define
$\La^{i+1}$ as the group of translations generated by all vertices
of $C^i$ (since we have chosen the origin in $A$, this makes
sense). Then, every vertex of $C^i$ is special in $C^{i+1}=(A,
W_{af}^{i+1})$, $W_{af}^{i+1}=\La^{i+1}\rtimes W$. By taking
$$
\La':= \bigcup_i \La^i, \quad W'_{af}:= \bigcup_i W^i_{af},
$$
we obtain the required special Coxeter complex $(A, W'_{af})$.
\qed

Given a point $x\in A=E^n$, we then let $W_{af,x}$ denote the
stabilizer of $x$ in $W_{af}$. Then the {\em type}, $type(x)$, of
the point $x$ is the $W$-conjugacy class of the {\em linear part}
of $W_{af,x}$ (i.e., the image under the canonical projection
$W_{af}\to W=W_{af}/A\cap W_{af}$). In particular, for special
Coxeter complexes, all vertices have the same type. Note that this
concept of type is weaker\footnote{Unless $W_{af}$ is transitive
on $A$.} than the standard (in the discrete case) notion of type
given by the coset $W_{af}\cdot x$ of $x$ in $A/W_{af}$
(equivalently, the intersection of the orbit $W_{af}\cdot x$ with
the fundamental alcove). However, in the case of nondiscrete
Coxeter groups, the latter concept does not seem to be of much use
since $A/W_{af}$ cannot be identified with an alcove or any
reasonable (commutative) geometric object.

In the spherical case, the notion of {\em type} is given by the projection
$$
\theta: S^{n-1}\to S^{n-1}/W=\De_{sph},
$$
where the quotient is the spherical Weyl chamber, a fundamental domain for $W\acts S^{n-1}$.

\medskip
We now assume that $n=2$, $A=E^2$ and $W=D_m$, the dihedral group
of order $2m$. Then there will be at most $div(m)+2$ types of
points in $(A, W_{af})$, where $div(m)$ is the number of divisors
of $m$. We choose the coordinates in $A$ so that the positive Weyl
chamber $\De$ of $W$ is bisected by the $x$-axis.

A {\em sub-wall} in $A$ is a codimension 1  cell in $A$, i.e., a
nondegenerate geodesic segment contained in a wall and bounded by two vertices.

\section{Metric concepts}\label{metric}

\begin{notation}
For a subset $Y$ in a metric space $X$ and $r\ge 0$, we let
$B_r(Y)$ denote the {\em open $r$-neighborhood} of $Y$ in $X$,
i.e.,
$$
B_r(Y):= \{x\in X: \exists y\in Y, d(x,y)< r\}.
$$
For instance, if $Y=\{y\}$ is a single point, then $B_r(Y)=B_r(y)$
is the open $r$-ball centered at $y$.
\end{notation}

A metric space $X$ is called {\em geodesic} if  every two points
in $X$ are connected by a (globally distance-minimizing) geodesic.
Most metric spaces considered in this paper will be geodesic; the
only occasional exceptions will be links of vertices in certain
metric cell complexes.

\begin{notation}
For a pair of points $x, y$ in a metric space $X$  we let
$\ol{xy}$ denote a closed geodesic segment (if it exists) in $X$
connecting $x$ and $y$. As, most of the time, we will deal with
spaces where every pair of points is connected by the unique
geodesic, this is a reasonable notation. We let $xy$ denote the
{\em open geodesic segment}, the {\em interior of} $\ol{xy}$:
$$
xy= \ol{xy}\setminus \{x, y\}.
$$

Geodesics will be always parameterized by their arc-length.
\end{notation}

{\bf $D$-isometries.} Let $X, X'$ be metric spaces and $D>0$ be a
real number. We say that a map $f: X\to X'$ is a {\em
$D$-isometry} if it satisfies the following:

1. $\forall x, y\in X, d(x,y)< D \Rightarrow d(f(x),
f(y))=d(x,y)$.

2. $\forall x, y\in X, d(x,y)\ge D \Rightarrow d(f(x), f(y))\ge
D$.

\medskip
The motivation for this definition comes from considering maps
where the diameter $D$ of the target is less than the diameter of
the domain, but we still would like to have maps $X\to X'$ which
are as close to being isometries as possible. We will use this
concept, however, even if $\diam(X')>D$.

The following properties of $D$-isometries are immediate:

\begin{itemize}

\item Each $D$-isometry is injective.

\item For each $D$-isometry $f$, its inverse $f^{-1}: f(X)\to X$
is also a $D$-isometry.

\item Composition of $D$-isometries is again a $D$-isometry.
\end{itemize}

{\bf CAT(0) and CAT(1) spaces.} We refer the reader for instance
to \cite{BH} for the detailed treatment of this material.

Let $X$ be a geodesic metric space. The space $X$ is called
CAT(0),  if the geodesic triangles in $X$ are {\em thinner} than
triangles in the Euclidean plane $E^2$. One can state this in
terms of quadruples of points in $X$. Let $x, y, z\in X$ (vertices
of a triangle $\Del(x,y,z)$). Let $m\in X$ be such that
$$
d(x,m)+d(m,y)=d(x,y)
$$
(i.e., $m$ will be a point on the side $\ol{xy}$ of the triangle).
Find points $\bar{x}, \bar{y}, \bar{z}$ in $E^2$ so that
$$
d(x,y)=d(\bar{x}, \bar{y}), \quad d(y,z)=d(\bar{y}, \bar{z}),
\quad d(z,x)=d(\bar{z}, \bar{x}),
$$
these are vertices of a {\em comparison triangle} $\Del(\bar{x},
\bar{y}, \bar{z})$ in $E^2$. Next, find a point $\bar{m}\in E^2$
so that
$$
d(x,m)=d(\bar{x}, \bar{m}), \quad d(m,y)=d(\bar{m}, \bar{y}),
$$
i.e., $\bar{m}$ belongs to the side of $\Del(\bar{x}, \bar{y},
\bar{z})$ which corresponds to the side $\ol{xy}$ of
$\Del(x,y,z)$. Then we require the CAT(0) inequality:
$$
d(z, m)\le d(\bar{z}, \bar{m}).
$$

\begin{figure}[tbh]
\centerline{\epsfxsize=4in \epsfbox{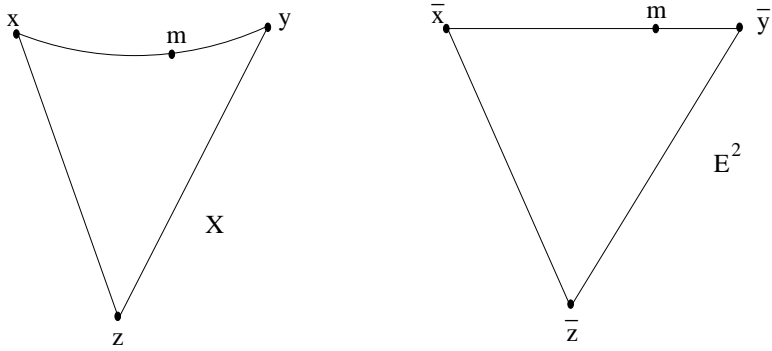}} \caption{\sl
Comparison triangle.} \label{f12}
\end{figure}

\medskip
Similarly, one defines CAT(1) spaces $Y$, except:

1.  Triangles of the perimeter $\le 2\pi$ in $Y$ are compared to
triangles on the unit 2-sphere.

2.  $Y$ need not be geodesic, however, if $p, q\in Y$ cannot be
connected by a geodesic, we require the distance between these
points be $\ge \pi$.

We will think of the distances in CAT(1) spaces as {\em angles}
and, in many cases, denote these distances $\angle(xy)$.

\medskip
The following characterization of 1-dimensional CAT(1) spaces will
be important:

A 1-dimensional metric space (a metric graph) is a CAT(1) space if
and only if the length of the shortest embedded circle in $X$ is
$\ge 2\pi$.

If $\Ga$ is a metric graph, where each edge is given the length
$\pi/m$, then the CAT(1) condition is equivalent to the assumption
that girth of $\Ga$ is $\ge 2m$.

\medskip
In a CAT(0) space, any pair of points is connected by a unique
geodesic, while in CAT(1) spaces this is true for points within
distance $<\pi$. A subset $Z$ in a CAT(0) space $X$ is called {\em
convex} if for every two points $p, q\in Z$ the geodesic
$\ol{pq}\subset X$ is contained in $Z$. A path-connected subset
$Z$ in a CAT(0) space $X$ is convex if and only if it is locally
convex, i.e., every point $z\in Z$ possesses a neighborhood in $Z$
which is convex in $X$.

\medskip
{\bf Spaces of direction.} Let $X$ be a CAT(0) space. The {\em
space of directions} $\Si_x(X)$ at a point  $x\in X$ is defined as
the set of germs of geodesics emanating from $x$. The distances in
$\Si_x(X)$ are given by the {\em angles} between the geodesics.

If $X$ is CAT(0), then $\Si_x(X)$ is CAT(1) for every $x\in X$.
Conversely, if $X$ is simply-connected, and $\Si_x(X)$ is CAT(1)
for every $x\in X$, then $X$ is CAT(0).

\begin{example}
Suppose that $X$ is a 2-dimensional Euclidean metric cell complex,
i.e., $X$ is obtained by gluing isometrically convex polygons in
$E^2$.  The space of directions $\Si_x(X)$ is just the link $\Ga$
of $x$ in $X$. To describe the metric on $\Si_x(X)$ note that a
small neighborhood of $x$ in $X$ is covered by flat Euclidean
sectors $S$  with the tip at $x$. The edges of $\Ga$ correspond to
maximal flat sectors $S$ (not containing any edges or vertices in
its interior). The angles of these maximal sectors at $x$ are the
lengths of the corresponding edges of $\Ga$.
\end{example}

\section{Geometry of spaces modeled on $(A,
W_{af})$}\label{spaces}

We refer to \cite{KL} for the detailed discussion.

{\bf Spaces modeled on Coxeter complexes.}

\begin{defn}\label{modeled}
A space  {\em modeled on the Coxeter complex $(A=E^n, W_{af})$} is
a metric space $X$ together with an atlas where charts are
isometric embeddings $A\to X$ and the transition maps are
restrictions of the elements of $W_{af}$. The maps $A\to X$ and
their images are called {\em apartments} in $X$. Note that (unlike
in the definition of an atlas in a manifold) we do not require the
apartments to be open in $X$. The same definition applies to
spaces modeled on the spherical Coxeter complex $(S^{n-1}, W)$,
only now the charts are maps of $S^{n-1}$.
\end{defn}

\begin{rem}
This definition generalizes to the case where we replace $W_{af}$
with an arbitrary group $G$ of affine transformations of $E^n$ and
assume that $A\subset E^n$ is a $G$--invariant affine subspace. In
this generality, one has to remove, of course, the assumption that
$X$ is a metric space and the charts are isometric embeddings.
Instead, one has to assume that for every two charts $\phi, \psi:
A\to X$,
$$
\psi^{-1} \phi(A)
$$
is a closed convex subset of $A$. Interesting examples of this
setup are provided by {\em tropical geometry}, where $G=\R^n
\rtimes SL(n,\Z)$.
\end{rem}

Therefore, all $W_{af}$-invariant (resp. $W$-invariant) notions defined in $A$, extend
to $X$. In particular, we will talk about vertices, walls, sub-walls, cells, etc. Note that the
intersection of finitely many cells in $X$ is again a cell.

We also will have a well-defined notion of {\em type} of a point
$x\in X$: It is given by computing the type of $x$ in an apartment
$A\subset X$ containing $x$. We retain the notation $type(x),
\theta(x)$ in the affine, resp. spherical, case.

\begin{ex}
If $W=D_m$ and $Y$ is modeled on $(S^1, W)$, then $Y$ has the
natural structure of a bipartite metric graph, where every edge
has the length $\pi/m$.

Accordingly, if $X$ is a CAT(0) space modeled on $(A, W_{af}=\La
\rtimes W)$, then for every $x\in X$, $Y=\Si_x(X)$ is a CAT(1)
metric space modeled on $(S^1, W')$, where $W'=type(x)$, a
subgroup of $W$.
\end{ex}

A space $X$ modeled on a Coxeter complex is called {\em locally
finite} if every point of $X$ has a neighborhood covered by
finitely many apartments.
A space $X$ modeled on $(A, W_{af})$ is said to be {\em finite} if it is covered by
finitely many apartments. (This, of course, does not mean finiteness of $X$ as a set or
finiteness of its diameter.)

\begin{notation}\label{complex}
We let $\c(A, W_{af})$ denote the collection of all finite spaces $X$ modeled on $(A, W_{af})$.
When  $(A, W_{af})$ is fixed, we will use the notation $\c$ for this collection. Similarly, we define
$\c_{lf}$ the collection of all locally finite spaces $X$ modeled on $(A, W_{af})$.
\end{notation}

A {\em finite subcomplex} $Y$ in $X\in \c_{lf}$ is a finite union
of cells in $X$. Note that such $Y$ need not be modeled on $(A,
W_{af})$.


\begin{convention}
In what follows, all spaces modeled on an affine Coxeter complexes are assumed to be CAT(0)
and all metric spaces modeled on spherical Coxeter complexes are assumed to be CAT(1).
\end{convention}

For a space $Y$ modeled on a spherical Coxeter complex $(S^{n-1},
W)$, two elements $\xi, \eta\in Y$ are called {\em antipodal} if
$\angle(\xi, \eta)=\pi$. Elements of $Y$ are {\em super-antipodal}
if $\angle(\xi, \eta)\ge \pi$. Accordingly, chambers $\si_1,
\si_2\subset Y$ are antipodal (rep. super-antipodal) if they
contain antipodal (resp. super-antipodal) regular elements. For
instance, if $Y$ consists of a single apartment, then chambers
$\si_1, \si_2\subset Y$ are antipodal if and only if the element
$w\in W$ sending $\si_1$ to $\si_2$ is the longest element of $W$,
i.e.,
$$
\si_2=-\si_1=\{-\xi: \xi\in \si_1\}.
$$

\medskip
Let $X$ be modeled on a 2-dimensional complex $(A, W_{af})$. A
{\em thick sub-wall} in $X$ is a sub-wall so that for every $x\in
\ga$, $\Si_x(X)$ is not a single apartment. In other words, $x$
has no neighborhood which is locally isometric to $A$. Of course,
not every sub-wall is thick.

Similarly, we say that $x\in X$ is a {\em thick vertex} of $X$ if
either $x$ locally separates $X$ or $x$ belongs to at least three
thick sub-walls $e_1, e_2, e_3$, so that
$$
e_i\cap e_j =\{x\}, i\ne j.
$$
Equivalently, the graph $\Si_x(X)$ is either disconnected or
contains at least 3 vertices of valence $\ge 3$.

\begin{ex}\label{examples}
1. Take two apartments $A_1, A_2$ and glue them at a vertex $v$.
Then $x$ is a thick vertex of the resulting space.

2. Take apartments $A_1, A_2, A_3, A_4$ and glue them along walls $L_1$, $L_2$
as follows:
$$
X_1:=A_1\cup_{L_1}\! A_2, X_2:= A_3 \cup_{L_2}\! A_4.
$$
Now, take an isomorphism $\phi: A_1\to A_2$, so that $\phi(L_1)$
crosses (transversally) $L_2$. Next, glue $X_1, X_2$ via $\phi$ to
form a space $X$. Then the vertex $v\in X$ corresponding to
$\phi(L_1)\cap L_2$ is a thick vertex in $X$.
\end{ex}

\medskip
A {\em Euclidean (or, affine) building} modeled on $(A, W_{af})$ is a CAT(0) space $X$ modeled on
$(A, W_{af})$ which satisfies the following two extra conditions:

Axiom 1. (``Connectedness'') Every two points $x_1, x_2\in X$ are contained in a common apartment.

\medskip
Axiom 2. (``Angle rigidity'') For every $x\in X$, the space of directions $Y=\Si_x(X)$ satisfies the following:
$$
\forall \xi, \eta\in Y, \quad \angle(\xi, \eta)\in W\cdot \angle(\theta(\xi), \theta(\eta)).
$$
Here $\theta: Y\to \De_{sph}$ is the {\em type projection}.

For instance, one can see that every (non-flat)
nonpositively curved symmetric space satisfies Axiom 1 but not Axiom 2.

A {\em spherical building} is defined similarly, except that it is
required to be modeled on a spherical Coxeter complex, be CAT(1)
and the Axiom 2 is not needed in this case (it follows from {\em
discreteness} of the Coxeter complex). It was proven in \cite{CL}
(and in an unpublished paper by B.~Kleiner)  that a
piecewise-Euclidean CAT(0) cell complex is a Euclidean building if
and only if all its links are spherical buildings.

\medskip
A building $X$ is called {\em thick} if every wall in $X$ is the intersection of
(at least) three half-apartments.

Below is an alternative set of axioms for Euclidean buildings  due
to Ann Parreau \cite{Parreau}. Replace Axioms 1 and 2 with:

Let $X$ be a CAT(0) space modeled on $(A, W_{af})$ 
so that:

\medskip
Axiom 1$^\prime$. Let $x_1, x_2\in X$ and $\si_i\subset
\Si_{x_i}(X)$ be spherical chambers, where we endow the spherical
buildings $\Si_{x_i}(X)$ with the $(S^1, type(x_i))$ structure.
Then there exists an apartment $A\subset X$ so that $x_i\in A,
\si_i\subset \Si_{x_i}(A)$, $i=1,2$.

\medskip
Axiom 2$^\prime$. For every $x\in X$ the following holds. Let
$\De_i\subset X, i=1,2$ be antipodal {\em  Weyl sectors} with the
tip $x$ (i.e., $\Si_{x}\De_i$ are antipodal in $\Si_x(X)$). Then
there exists an apartment $A\subset X$ containing $\De_1, \De_2$.

\medskip
Here a {\em  Weyl sector} in $X$ is a subset $\De'$ of an
apartment $A\subset X$, so that
$\De'$ is a parallel translate of a Weyl chamber in $A$. Note that if
the tip $x$ of $\De'$ is a special vertex,
then $\De'$ is a Weyl chamber in $X$.

\begin{rem}
Parreau's numbering of axioms is different from ours. She also does not require the CAT(0)
condition, but has other axioms requiring intersections of pairs of apartments to be closed
and convex. These conditions follow from the CAT(0) property.
\end{rem}

 We will be using Axioms 1$^\prime$ and 2$^\prime$
since verifying angle rigidity would require proving a slightly
weaker form of Axiom 1$^\prime$ anyhow. We refer the reader to
Tits' original paper \cite{Tits3} and to Appendix 3 of Ronan's
book \cite{Ronan2} for alternative axiomatization of nondiscrete
Euclidean buildings.

\medskip
{\bf Morphisms.}

A {\em morphism} between two spaces $X, X'$ modeled on $(A,G)$,
where $A=E^n, G=W_{af}$ or $A=S^{n-1}, G=W$, is defined as a map
$f: X\to X'$  which, being written  in ``local coordinates given
by apartments'' is the restriction of elements of $G$.

More generally, given two spaces $X, X'$ modeled on $(A,G)$, and
subsets $Y\subset X, Y'\subset X'$, a {\em morphism} $Y\to Y'$ is
a map $f: Y\to Y'$, which, in local coordinates, appears as the
restriction of an element of $G$. In other words, for every point
$y\in Y, y'=f(y)$, and for every pair of apartments $\phi: A\to X,
\psi: A\to X'$ (with images containing $y, y'$ respectively) the
(partially defined) map
$$
\psi^{-1}\circ f\circ \phi:  \phi^{-1}(Y)\subset A \to  \psi(Y') \subset A
$$
is the restriction of an element of $G$.

It follows that, if a geodesic $\ga\subset X$ is contained  in
$Y$, then its image under a morphism $Y\to Y'$ is a
piecewise-geodesic path. Therefore, every morphism $Y\to Y'$
induces a morphism (again denoted $f$ to simplify the notation)
$$
f: \Si_y Y \to \Si_{y'} Y', \quad y\in Y, \quad y'=f(y)\in Y'.
$$

\medskip
{\bf Weak isometries.}  Let $Z$ be a subcomplex in a metric space
$U$ modeled on a spherical Coxeter complex $(S^{n-1}, W)$. We
always endow $Z$ with the restriction of the path-metric on $Z$
(not the induced path-metric!). Let $V$ denote the vertex set of
the complex $(S^{n-1}, W)$. Set
$$
D:= \max \{d(x, w(x)): x\in V, w\in W\}.
$$
Suppose now that $n=2$ and $W=D_m$. We set $\del=\pi -\pi/m$.
Then, clearly,
$$D=
\begin{cases}
\pi         & \text{if $m$ is even}\\
\pi-\pi/m=\del & \text{if $m$ is odd} \\
\end{cases}
$$

More generally, for arbitrary finite Coxeter groups $W$ we have:
$D\le \pi$ with strict  inequality if and only if the Coxeter
number of $W$ is odd. This is the case when $W$ is a product of
symmetric groups $S_m$ with odd $m$ and dihedral groups $D_m$ with
odd $m$. For instance, if $W=S_m$ with odd $m$, then
$D=\arccos(-(m-1)/(m+1))$.

\medskip
We then say that a {\em weak isometry} between subcomplexes
$Z\subset U, Z'\subset U'$ is a surjective morphism $Z\to Z'$
which is a $D$-isometry. It is easy to see that in the case
$W=D_m$ and $n=2$, the latter condition is equivalent to being a
$\del$-isometry.

\medskip
We now define weak isometries between subsets of spaces modeled on
affine Coxeter complexes. Let $X, X'$ be modeled on the affine
Coxeter complex $(A, W_{af})$.

1. We first consider the easier case when $Y\subset X, Y'\subset
X'$ are subcomplexes. Then a {\em weak isometry} $Y\to Y'$ is a
morphism $Y\to Y'$ which sends cells isometrically to cells  and
induces a weak isometry
$$
\Si_y(Y)\to \Si_{y'}(Y')
$$
for every $y\in Y, y'=f(y)\in Y'$. This is the concept akin to the
notion of an isometric embedding of Riemannian manifolds (rather
than isometry between metric spaces).

2. Now, suppose that $Y\subset X, Y'\subset X'$ are {\em locally
conical subsets}, i.e., each point of $Y, Y'$ admits a
neighborhood in $Y, Y'$ which is a subcomplex in $X, X'$.
Therefore,  for every $y\in Y, y'\in Y'$,
$$
\Si_y(Y)\subset \Si_y(X), \quad \Si_{y'}(Y')\subset \Si_{y'}(X')
$$
are subcomplexes.

Then a {\em weak isometry} $Y\to Y'$ is a surjective morphism $f:
Y\to Y'$ inducing a weak isometry
$$
f: \Si_y(Y)\to \Si_{f(y)}(Y'),
$$
for every $y\in Y$.

We will use the name {\em weakly isometric embedding} for the map $f: Y\to X'$,
when we do not want to specify the image of $f$.

\medskip
{\bf Local structure of a locally finite space $X$ modeled on $(A,
W_{af})$.}

\begin{defn}
(Conical sets, cf. \cite{CL}.)  Let $Z\subset X$ be a subset and $o\in Z$. Then
we say that $Z$ is {\em conical}  with respect to $o$ if:

1. For every $x\in U$, $\ol{xo}\subset Z$ (i.e., $Z$ is star-like
with respect to $o$).

2. Geodesics starting at $o$ do not branch in $Z$: If $\ga_i=\ol{o
x_i}$, $i=1, 2$ are geodesics of the same length in $C$, and the
germs of $\ga_1, \ga_2$ agree near $o$, then $x_1=x_2$.
\end{defn}

Let $Y$ be a metric graph. We define the {\em cone over $Y$}
$Cone(Y)$ as follows. For each edge $e$ of $Y$ of length $\al$,
let $Cone(e)$ be the (infinite) sector in $\R^2$ with the tip $o$
(at the origin) and the angle $\al$. If $v$ is a vertex of $e$,
the cone $Cone(v)$ the ray bounding $Cone(e)$. We now glue all the
sectors $Cone(e)$ via isometries
$$
Cone(v)\subset Cone(e)\to  Cone(v)\subset Cone(e')
$$
whenever the edges $e, e'\subset Y$ share the vertex $v$. The
result is $Cone(Y)$. There is a more general definition of a
metric cone over a CAT(1) space $Y$, obtained by taking an
appropriate metric on $Y\times [0,\infty)/Y\times \{0\}$, see
\cite{KL,N}, but we will not need it.

\medskip
Now, if $Z$ is a conical subset in $X$ (modeled on $(E^2, W_{af})$),
then at every point $o\in Z$,
the set $Z$ is isometric to a star-like subset of the
Euclidean cone
$$
Cone(Y), \quad Y=\Si_o(X).
$$

Since $X$ is locally finite, each point $x\in X$ has a convex
neighborhood which is conical with respect to $x$. In particular,
every subcomplex in $X$ is locally conical.

The following lemma describing the geometry of $X$ near sub-walls
generalizes the above observation. We first need

\begin{defn}
Let $X_0$ be a building which contains at most one thick wall $L$,
i.e., $X_0$ is obtained by attaching several half-apartments along
$L$. Then, for each closed subsegment $\ga\subset L$ and $\eps>0$
we consider $B_\eps(\ga)\subset X_0$. We will refer to
$B_\eps(\ga)$ as an {\em open book with the binding $\ga$}.
\end{defn}

\begin{figure}[tbh]
\centerline{\epsfxsize=4in \epsfbox{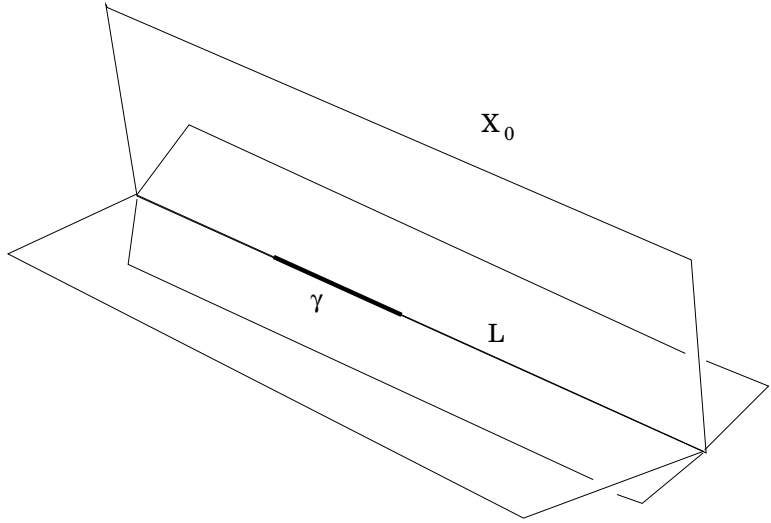}}
\caption{\sl Open book with the binding $L$.}
\label{f5}
\end{figure}

\begin{lem}
Let $X$ be locally finite, modeled on $(A, W_{af})$. Let
$\ga\subset X$ be a sub-wall which does not cross (transversally)
any thick sub-walls. (We allow $\ga$ to be a sub-wall.) Then
there exists $\eps>0$ so that $B_\eps(\ga)$ is conical with
respect to every $o\in \ga$. Moreover, $B_\eps(\ga)$ is isomorphic
to an {\em open book with the binding $\ga$}.
\end{lem}
\proof 1. Since $X$ is CAT(0), the neighborhood $C=B_\eps(\ga)$ is
convex. In particular, it is star-like with respect to every point
$o\in \ga$.  If two geodesics in $C$ ``branch'', then branching
occurs at a point of $\ga$, since elsewhere, $C$ is locally
isometric to $A=\R^2$. If a geodesic $\be=\ol{ox}\subset C$
intersects $\ga$ at a point different from $o$, then it is
entirely contained in $\ga$. The non-branching of geodesics
follows.

2. Pick a point $o\in \ga$ which is not an end-point. We identify
$X_0$ with $Cone(\Si_o(X))$. We have the map
$$
h: C\to X_0
$$
which sends $x\in C\setminus \{o\}$ to the pair $(\xi, r)$, where
$\xi$ is the direction of $\ol{ox}$ in $\Si_o(X)$ and $r=d(o,x)$.
We leave it to the reader to verify that $h$ is an isometry of $C$
to an open book in $X_0$. \qed

\medskip
{\bf Unit tangent bundles.} We observe that every $X$ modeled on
$(A, W_{af})$ has extendible geodesics: If $\ga: [0, T]\to X$ is a
geodesic, it extends to a geodesic ray $\ga: \R_+\to X$. This
follows since every  $\xi\in \Si_x(X)$ has an antipodal direction
(coming from an apartment in $\Si_x(X)$ containing $\xi$).

For a   subcomplex $P\subset X$ (or, more generally, a locally
conical subset), we define the ``unit tangent bundle'' $\Si(P)$ as
the space of pair $(x, \xi), x\in X, \xi\in \Si_x(P)\subset
\Si_x(X)$. Then $\Si(P)\subset \Si(X)$. We topologize $\Si(X)$ as
follows. Let $Ray(X)$ denote the space of geodesic rays in $X$
with the topology of uniform convergence on compacts. Then
$\Si(X)$ is the quotient of $Ray(X)$, where we identify rays
emanating from $x$ having the same germs at $x$. We give $\Si(X)$
the quotient topology and $\Si(P)$ the subspace topology. The
following are elementary properties of $\Si(X)$:

\begin{itemize}
\item $\Si(X)$ is Hausdorff.

\item The projection $\pi: \Si(X)\to X$ is continuous.

\item If $P\subset X$ is closed then $\Si(P)\subset \Si(X)$ is closed. (Note that in general
$\Si(P)\ne \pi^{-1}(P)$.)

\item If $X$ is locally finite then the map $\pi$ is proper
(Arcela-Ascoli theorem). In particular, $\Si(X)$ is locally
compact in this case.
\end{itemize}

\medskip
{\bf Developing maps.} Let $U\subset X$ (which is modeled on
2-dimensional $(A, W_{af})$) be a connected subset which is
homeomorphic to an open 2-dimensional manifold. We also assume
that each point $x\in U$ has a neighborhood which is contained in
a single apartment. We thus obtain an {\em open} atlas $\{(U_\al,
\phi_\al)\}$ on $U$, whose charts
$$
\phi_\al: U_\al\to \phi_\al(U_\al)\subset A
$$
are  isomorphisms. Thus, the transition maps for this atlas are
restrictions of the elements of $W_{af}$. This gives us a
geometric structure on $U$ modeled on $(A, W_{af})$ and,
therefore, we obtain a multi-valued developing map $f: U\to A$
which is, locally, a morphism (see e.g. \cite{G}).

\begin{lem}
Suppose that $U$ is convex in $X$. Then $f$ is single-valued and
injective.
\end{lem}
\proof Since $U$ is convex, it is simply-connected. Hence, $f$ is
single-valued. In order to check injectivity, note that, being a
local isomorphism, $f$ is a local isometry. Let $x, y\in U$ be
distinct points. Since $U$ is convex, the (unique) geodesic
$\ga=\ol{xy}\subset X$ is contained in $U$. Since $f$ is a local
isomorphism, it sends local geodesics to local geodesics. Since a
local geodesic in $A$ is a global geodesic, $f(\ga)\subset A$ is a
straight-line segment with the end-points $f(x), f(y)$ and the
length equal to that one of $\ga$. Therefore, $f(x)\ne f(y)$. \qed

\begin{defn}
A subset $P\subset X$ is called {\em planar} if there exists a
weakly  isometric embedding $f: C\to A$.
\end{defn}

\section{Corridors and tiled corridors}\label{corridors}

In this section we will assume that the Coxeter complex $(A,
W_{af})$ is {\em special}, the group of translations of $W_{af}$
is dense in $\R^2$, and spaces $X$ modeled on this complex are
locally finite. Below we define {\em tiles} and {\em corridors},
which are the main tools for proving a combinatorial convexity
theorem \ref{convex}, which, in turn, is the key for constructing
the building $X$ as in Theorem \ref{main}.

\medskip
A {\em model tile} in a Coxeter complex $(A,W_{af})$ is a
parallelogram $D$ whose edges are sub-walls parallel to the walls
bounding $\De$.  In particular, the angles of $D$ are of the form
$\pi/m, \pi-\pi/m$. We will regard a single sub-wall, with the
slope $\pm \tan(\pi/2m)$, as a (degenerate) tile. We let
$v_-=v_-(D), v_+=v_+(D)$ denote the left and right vertices of
$D$. Then $D$ is uniquely determined by $v_\pm$.

We now define {\em model corridors}  in $A$. We will frequently
omit the adjective ``model'' to simplify the terminology. Every
model corridor $C$ is determined a pair of piecewise-linear
functions $u_\pm: [a, b]\subset \R \to \R$ satisfying the
following:

1. $u_+(a)=u_-(a), u_+(b)=u_-(b)$. The points $v_-=(a,u_\pm(a)),
v_+=(b,u_\pm(b))$ are (distinct) vertices of $A$: They will be
called the {\em extreme vertices} (or {\em extremes}) of the
corridor $C$.

2. $u_-(x)\le u_+(x)$ for all $x\in [a,b]$.

3. The slopes of the graphs of $u_\pm$ are of the form $\pm
\tan(\pi/2m)$.

The function $u_+$ is called {\em upper roof function} and $u_-$
is called {\em lower roof function.} The slope of $u_+$ at $v_-$
is $\tan(\pi/2m)$ and at $v_+$ it is $-\tan(\pi/2m)$, for $u_-$ it
is the other way around.

\begin{rem}
In Section \ref{generalized}, we will generalize the above
concepts to allow tiles of infinite diameter and roof functions
defined on infinite intervals.
\end{rem}

\begin{defn}
A {\em (model) corridor $C$} is the closed region between the
graphs of the functions $u_+, u_-$. A corridor $C$ is called {\em
simple} if either $u_-(x)<u_+(x)$ on the open interval $(a,b)$ or
$u_+\equiv u_-$ on $[a,b]$. In other words, simple corridors are
the ones which are homeomorphic to a closed disk or isometric to a
segment.

A corridor is {\em degenerate} if it is a sub-wall, i.e., is a
single segment. Thus, every degenerate corridor is simple.
\end{defn}

\begin{figure}[tbh]
\centerline{\epsfxsize=4in \epsfbox{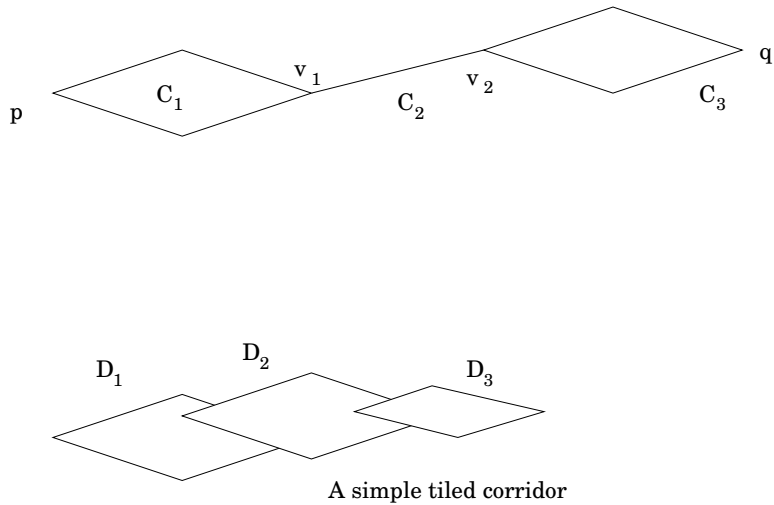}} \caption{\sl
Corridors.} \label{f11}
\end{figure}

Figure \ref{f11} provides examples of corridors.

We note that each corridor is a union of finitely many simple
corridors, every two of which intersect in at most one point, an
extreme vertex of both corridors. These intersection points will
be called {\em cut-vertices} of $C$, as they separate the corridor
$C$. Away from these points, a corridor $C$ is a manifold with
boundary (of dimension 1 or 2). This decomposition of $C$ into
maximal simple sub-corridors $C_i$ is unique; the sub-corridors
$C_i$ are called the {\em components} of $C$.

We let $\D_+ C, \D_- C$ denote the graphs of $u_+, u_-$. Then $\D
C=\D_+ C\cup \D_- C$.

\medskip
We note that each corridor with the extremes $v_-, v_+$ is
contained in the tile with the vertices $v_-, v_+$. It follows
immediately from the definition that each edge of the polygon $C$
is a sub-wall and each vertex is a vertex of $(A,W_{af})$.

Despite the fact that a corridor is not convex (unless it is a
single tile), it satisfies the following weak convexity property:

Let $L$ be a line in $A$ with the slope $\pm \tan(\pi/2m)$. Then
$L\cap C$ is either empty or is a segment.

\medskip
Although we will not need this,  one can characterize simple
nondegenerate  corridors in $A$ geometrically (without using the
roof functions) as follows:

\begin{enumerate}
\item $C$ is a subcomplex in $(A, W_{af})$.

\item $C$  has only finitely many edges and is homeomorphic to the
closed 2-disk.

\item All edges of $C$ have slopes $\pm \tan(\pi/2m)$.

\item $C$ does not have any two consecutive (interior) angles each
$> \pi$.
\end{enumerate}

The following lemma is elementary and is left to the reader.

\begin{lem}\label{approx}
Suppose we are given $\eps>0$ and a continuous function $u:
[a_1,a_2]\to \R$, so that $(a_i, u(a_i))$ are vertices of the
Coxeter complex $(A, W_{af})$. Then there are upper and lower roof
functions $u_\pm: [a_1,a_2]\to \R$ so that:

1. $u_\pm(a_i)=u(a_i), i=1, 2$.

2. $\|u- u_\pm\|<\eps$.
\end{lem}

\begin{lem}\label{limit}
Let $C_s$ ($s\in \R_+$) be an increasing family of corridors with
common extreme points $v_-, v_+$: $C_s\subset C_t \iff s\le t$.
Assume also that the numbers of edges in the corridors $C_s$ are
uniformly bounded from above. Then
$$
C=\ol{\bigcup_{s} C_s}
$$
is a corridor in $(A, \ol{W_{af}})$ with the extreme points $v_-,
v_+$.
\end{lem}
\proof Let $u_{\pm,s}$ be the upper/lower roof function for the
corridor $C_s$. Since the corridors $C_s$ are increasing, we
obtain:
$$
u_{+,s}\le u_{+,t}, \quad u_{-,s}\ge u_{-,t}, \quad \forall s\le
t.
$$
Moreover, all the roof functions in question are uniformly bounded
(from above and below) by the roof functions of the tile with the
vertices $v_-, v_+$. It is then immediate that there exist limits
$$
u_\pm= \lim_{s\to\infty} u_{\pm,s}
$$
and these limits are roof functions with respect to $(A,
\ol{W_{af}})$. \qed

\medskip
It is clear from the definition that every point in a corridor $C$
belongs to the interior (relative to $C$) of a tile $D\subset C$.
By  compactness, we can therefore cover the corridor $C$ by
finitely many tiles.

\begin{defn}
A {\em tiled corridor} in $A$ is a finite set of tiles
$N=\{D_1,...,D_k\}$ whose union $|N|$ is a corridor $C$ in $A$, so
that $N$ is minimal with this property, i.e., any proper subset
$N'\subset N$ will have strictly smaller union. When convenient,
we will abuse the notation and conflate $N$ and $C=|N|$. We set
$\# N:=k$.
\end{defn}

\medskip We now extend our discussion to corridors in $X$, a space
modeled on $(A,W_{af})$.

\begin{defn}
A {\em tile} in $X$ is a subcomplex isomorphic to a (model) tile
in $A$.
\end{defn}

We now define corridors in $X$. They will satisfy two axioms C1
and C2, the first of which is easy to state,  the second axiom
will take more work.

\medskip
{\bf Axiom C1.} Every corridor in $X$ is a subcomplex $C\subset X$
weakly isometric to a model corridor in $A$.

\medskip
In particular, if $v$ is a cut-vertex of $C$, it disconnects $C$
into two sub-corridors $C_1, C_2$, so that $\Si_v(C_i)$ is a
chamber or a vertex in $\Si_v(X)$, and the distance between these
is at least $\pi- \pi/m$.

Using the embedding $f: C\to A$, we will be identifying corridors
$C$ in $X$ with corridors in $A$. Suppose $C\subset X$ is a
corridor with extreme vertices $v_-, v_+$.  Since $f(C)$ is
contained in a tile with the extreme vertices $f(v_-), f(v_+)$, we
obtain

\begin{lem}\label{diam}
$C\subset B_d(\ol{v_- v_+})$, where $d=d(v_-, v_+)$.
\end{lem}

\begin{figure}[tbh]
\centerline{\epsfxsize=4in \epsfbox{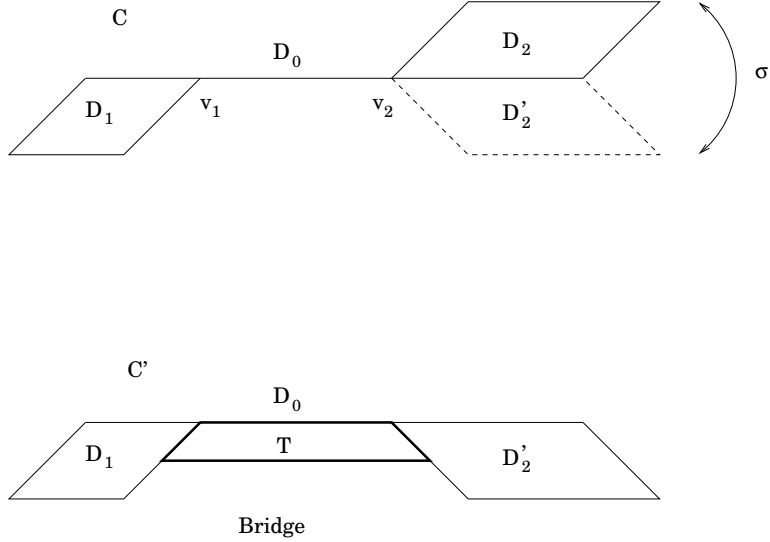}} \caption{\sl
Corridor, ``twisted corridor'' and a bridge. Here $D_0=\ol{v_1
v_2}$.} \label{f8}
\end{figure}

We now introduce the second corridor axiom. To motivate the
definition, consider a model corridor $C=D_0 \cup D_1\cup
D_2\subset A$ as in Figure \ref{f8}, where $D_0$ is a degenerate
tile. Let $\si\in W_{af}$ be the reflection in the wall containing
$D_0$. Then the union
$$
C'=D_0 \cup D_1\cup D_2' \subset A, \quad \hbox{where~~~} D_2'=\si(D_2),
$$
is not a model corridor in $A$, but satisfies Axiom C1: The
isomorphism $f: C'\to C$ is given by the map which fixes $D_1\cup
D_0$ pointwise and sends $D_2'$ to $D_2$ via $\si$. For a variety
of reasons, we would like to exclude $C'$ (and other similar
``twisted corridors'') from being corridors in spaces $X$ modeled
on $(A, W_{af})$. The most obvious reason (although not the
critical one) is that we would like every corridor in $A$ to be
the image of a model corridor under some $w\in W_{af}$. The key
reason for excluding such $C'$ will become clear in the proof of
Lemma \ref{wmax}.

One can tell apart a ``twisted corridor'' $C'\subset A$ from a
model corridor by observing that there exists a {\em trapezoid}
$T$ {\em connecting} the tiles $D_1, D_2'$ as in Figure \ref{f8},
while a model corridor cannot admit such trapezoids since its
existence would mean that boundary slopes of a model corridor
would take at least 3 different values. Equivalently, there is a
nontrivial continuous family of segments parallel to $D_0$ with
end-points contained in $C'$. (These segments foliate the
trapezoid $T$.)

\begin{defn}
A {\em trapezoid} in $A$ is a nondegenerate quadrilateral bounded
by sub-walls and having the angles $\pi/m, \pi/m, \pi-\pi/m,
\pi-\pi/m$.
\end{defn}

 We will say that the set of polygons
$D_1, D_0, D_2', T$ is a {\em bridge} in $A$. By abusing the notation, we will refer to the
union of polygons in a bridge as a bridge as well.

We define a {\em bridge} in a space $X$ (modeled on $(A, W_{af})$) as a set of polygons whose union is
isometric to a bridge in $A$.

We are now ready to state the second corridor axiom:

{\bf Axiom C2}. A corridor in $C$ cannot admit a ``bridge'', i.e.,
$C$ cannot contain a degenerate sub-corridor $D_0$ which, together
with tiles $D_1, D_2'\subset C$, and a trapezoid $T\subset X$ form
a bridge. Equivalently, if $D_0$ is a degenerate irreducible
component of $C$, there is no nontrivial family of segments in $X$
parallel to $D_0$ whose end-points are contained in $C$.

\medskip
Given a corridor $C\subset X$ and its weak isometry $f: C\to
C'=f(C)\subset A$, we declare the {\em boundary} $\D C$ of $C$, to
be the preimage $f^{-1}(\D C')$. It is clear that the boundary is
independent of the embedding $f$. We then similarly extend to
corridors in $X$ all the concepts defined for model corridors.

\begin{lem}\label{isomembed}
Suppose that $C\subset X$ is a corridor and a morphism $\iota:
X\to X'$ is an isometric embedding. Then $\iota(C)$ is again a
corridor in $X'$.
\end{lem}
\proof Let $f: C\to C'\subset A$ be a weakly isometric embedding. Then
$\iota\circ f^{-1}: C'\to \iota(C)$
is again a weakly isometric embedding. This verifies Axiom C1.

Suppose that $\iota(C)$ admits a bridge $D_0\cup D_1\cup D_2 \cup T$. Then $T$ is foliated by
parallel geodesic segments whose endpoints belong to $D_1\cup D_2\subset \iota(X)$. Since $\iota(X)$ is
convex in $X'$, it follows that each of these geodesics is contained in $\iota(X)$.
In particular, $C\subset X$ admits a bridge. Contradiction. \qed

\begin{rem}\label{nonconvex}
It important to note that, unless $X$ is a building, if $C\subset
X$ is a convex corridor, its image $f(X)\subset A$ is not
necessarily convex. This follows from the fact that a weakly
isometric embedding is not in general an isometric embedding. The
reader can construct examples of convex corridors with nonconvex
images $f(C)$ using the example of a space $X$ described in
Example \ref{examples}, Part 2.
\end{rem}

\section{Two orders on tiled corridors}\label{orders}

In this section, $X\in \c_{lf}$ (as in Notation \ref{complex}) is
a space modeled on $(A, W_{af})$.

\begin{defn}
Let $N=\{D_i, i\in I\}$ and $N'=\{D_j', j\in J\}$ be tiled
corridors in $X$ with $C=|N|, C'=|N'|$. We say that $N\ll N'$ if:

1. $C$ and $C'$ have the same extreme vertices.

2. For every $D_i\in N$  there exists $D_j'\in N'$ so that
$D_i\subset D_j'$.

3. $f_C=f_{C'}|_{C}$, where $f_C, f_{C'}$ are weak isometries to
model corridors in $A$.

 We say that $N<N'$ if $N\ll N'$ and the cardinality of $J$ is at
 most the cardinality of $I$, i.e., $\# N'\le \# N$.
\end{defn}

It is easy to check that both $<$ and $\ll$ are (partial) orders.

\begin{rem}
On can replace the order $\ll$ with the inclusion relation between
the corridors $C$, but the proofs will become somewhat  more
complicated.
\end{rem}

A tiled corridor $N$ in $X$ is {\em weakly maximal} if it is
maximal with respect to the order $<$. A tiled corridor $N$ in $X$
will be called {\em maximal}  if it is maximal with respect to the
order $\ll$.

\medskip
Let $X$ be modeled on $(A, W_{af})$. We will use the notation
$\ol{X}$ for the same metric space $X$, but regarded as a space
modeled on $(A, \ol{W_{af}})$, where $ \ol{W_{af}}=A\rtimes W$ is
the closure of $W_{af}$. Then

\begin{prop}\label{wmax}
For every tiled corridor $N_0$ in $X$ there exists a tiled
corridor $N$ such that $N$ is weakly maximal in $X$ and $N_0<N$.
\end{prop}
\proof Let $N_0=\{D_{0,0},...,D_{k,0}\}$.

1. We will first prove existence of a weakly maximal corridor in
$\ol{X}$; in the second part of the proof we will verify that this
corridor is actually a corridor in $X$. By Zorn's lemma, it
suffices to consider a maximal totally ordered (with respect to
the order $<$) family of tiled corridors
$$
\{N_s, s\in J\subset \R_+\}, $$ and to show that this family
contains a maximal element. Without loss of generality we can
assume all tiled corridors have the same cardinality. Then each
tiled  corridor $N_s$ is of the form
$$
\{D_{0,s},...,D_{k,s}\},
$$
where $D_{i,s}\subset D_{i,t}, s\le t$.

Set $C_s:= |N_s|$ and let $f_s: C_s\to A$ be the isomorphism onto
a model corridor $C_s'\subset A$, so that
$$
f_t|_{C_s}= f_s, \quad s\le t.
$$
Then we obtain an increasing family of (model) tiled corridors
$N_s'=f(N_s)$.

By Lemma \ref{limit}, the closure $C'$ of the union $\cup_s C'_s$
is again a corridor in $A$. Defining a tiled corridor
$$
N':= \{D'_{0},...,D'_{k}\}, \quad \hbox{where}\quad  D'_{i}:=
\ol{\cup_s D'_{i,s}}, \quad i=0,...,k.
$$
we get $C'=|N'|$. Similarly, set
$$
N:= \{D_{0},...,D_{k}\}, \quad \hbox{where~~} D_{i}:= \ol{\bigcup_s D_{i,s}},
\quad i=0,...,k, \quad C:=|N|.
$$

\begin{rem}
We may have to omit some of the tiles in $N$ and $N'$ in order to
obtain minimal tiled corridors $N, N'$ with $|N|=C, |N'|=C'$. In
order to simplify the notation we will ignore this issue.
\end{rem}

Then the maps $f_s$ yield a map $f: C\to C', f(D_i)=D_i',
i=0,...,k$. We claim that $C$ is a corridor in $\ol{X}$.

We first verify that $f$ is a weak isometry. It is clear that the
restriction of $f$ to every simple sub-corridor is an
isometry\footnote{with respect to the induced path-metric} to its
image, as the limit of isometries, is an isometry. In particular,
the restriction of $f$ to each simple sub-corridor is a weak
isometry.

Therefore, we have to consider the behavior of $f$ at the
cut-vertices. If $v$ is a cut-vertex of $C$ then it is a cut-point
of every $C_{s}, s\in \R_+$. If a cut-vertex $v$ of $C$ was
separating nondegenerate sub-corridors  for some $C_{s}$, the same
is true for all $t\ge s$. Moreover, in this case
$$
\Si_{v}(C_t)= \Si_{v}(C),\quad \Si_{v'}(C'_t)= \Si_{v'}(C'), \quad
v'=f(v), t\ge s,
$$
and
$$
f_t: \Si_{v}(C_t) \to \Si_{v'}(C'_t)
$$
is independent of $t$ for $t\ge s$. Thus, $f$ is a weak isometry
at $v$ in this case.

Suppose now that $v$ belongs to a degenerate component, say,
$D_{0,s}$,  of $C_s$ for all $s$. Note that the point $v$ need not
be an end-point  of this geodesic segment.

\begin{rem}
This is the part of the proof where the corridor Axiom C2  will be
used.
\end{rem}

We will consider only the case when $\Si_v( C_s)$ is {\em not} of the form
$$
\{\xi_1, \xi_2\}, \quad \hbox{where~~} \angle(\xi_1, \xi_2)=\del=\pi- \pi/m
$$
for sufficiently large $s$. The remaining case will be similar and is left to the reader.

If $v$ is not a vertex of $D_{0,s}$
we let $\xi_i\in \Si_v X$ ($i=1,2$) denote the antipodal directions tangent to $D_{0,s}$.
Otherwise, we define $\xi_i\in \Si_v X$ so that one of them (say, $\xi_1$) is
tangent to $D_{0,s}$ and the other is a fixed direction which
belongs to $\Si_v(C_s)$ (for all $s$) and is super-antipodal to $\xi_1$:
$$
\angle(\xi_1, \xi_2)\ge \pi.
$$
By the construction of $C$, the germ $\Si_{v}(C)$ at $v$ consists
of two components $\si_{i}$ containing $\xi_i$ ($i=1,2$). Each
$\si_i$ is either a chamber or equals $\xi_i$. Therefore, the
distance between $\si_1, \si_2$ in $\Si_v X$ is at least $\pi -
\frac{2\pi}{m}$. If the distance is $\ge \del=\pi - \pi/m$ then
$f$ is indeed a weak isometry at $v$ and we are done. Otherwise,
the distance equals   $\pi - \frac{2\pi}{m}$ and both $\si_i$ are
chambers $\ol{\xi_i \eta_i}$: $\si_i= \Si_v (D_{i}), i=1, 2$.
Moreover, we necessarily have $\angle(\xi_1, \xi_2)=\pi$. In
particular, $v$ separates nondegenerate sub-corridors in $C$.

\begin{figure}[tbh]
\centerline{\epsfxsize=4in \epsfbox{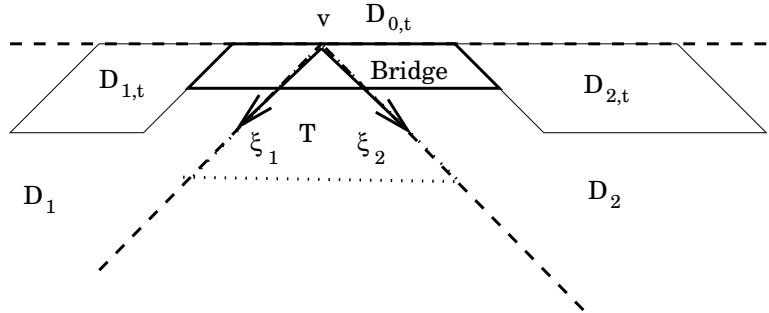}}
\caption{\sl Building a bridge between $D_{1,t}$ and $D_{2,t}$.}
\label{f9}
\end{figure}

Since the distance between $\si_1, \si_2$ is $\pi - \frac{2\pi}{m}$, there exists a flat triangle
$T\subset X$ with the tip $v$ and the angle $\pi - \frac{2\pi}{m}$ at $v$, so that
$$
\Si_v(T)= \{ \eta_1, \eta_2\}
$$
Combining the germ of this flat triangle with the germs (at $v$)
of $D_{i}, i=1,2$, we obtain a bridge in $X$ containing $D_{i,t},
i=0, 1, 2$ for sufficiently large $t$. This contradicts Axiom C2
for the corridors $C_t$. See Figure \ref{f9}.

The same argument excludes existence of a bridge for $C$: one can regard the triangle $T$
in the above argument as a degenerate trapezoid.

Since the family of tiled corridors $\{N_s\}$ was chosen to be
maximal, it shows that $N=N_s$ for all sufficiently large $s$ and
$N$ is a weakly maximal tiled corridor in $\ol{X}$.

\medskip
2. Note that even if all the corridors $C_s$ were corridors in
$X$, then $C$ might be a corridor only in $\ol{X}$, as the walls
of $(A, W_{af})$ are dense among the walls of $(A, \ol{W_{af}})$.
Suppose that one of the edges $e$ of $C$ is not contained in a
wall of $X$. Since its image under $f$ has the slope $\pm
\tan(\pi/2m)$, it follows that $e$ cannot contain any vertices of
$X$.

\begin{rem}
This is the only place of the paper where we {\em truly need} the
fact that $(A, W_{af})$ is special: All other arguments could be
modified to handle the general case.
\end{rem}

In particular, $e$ contains no extreme points of $C$ and no thick
vertices of $X$. Therefore, all the tiles of $N$ which overlap
with the edge $e$ can be slightly expanded in the direction of
this edge to form a new corridor $\hat{C}=|\hat{N}|$ which is
strictly larger than $C$ (it will have an edge parallel to $e$).
Thus, $N<\hat{N}, N\ne \hat{N}$ contradicting maximality of $N$.
\qed

\medskip
Our next goal is to prove existence of maximal tiled corridors
with respect to the order $\ll$. In view of the previous
proposition, it suffices to restrict our discussion to weakly
maximal tiled corridors.

Let $D, D'\subset A$ be model tiles; we denote the (closed) edges
of these tiles $e_1,...,e_4$ and $e_1',...,e_4'$ respectively, so
that the edges $e_i, e_i'$ correspond to each other under the
dilation
 mapping $D$ to $D'$. Given a subset $E\subset A$, we define
 $$
I=I(D,D', E):=\# \{ i\in \{1,...,4\}| \exists a\in E\cap e_i \cap
e_i'\},
 $$
which is {\em the number of sides of the tiles $D, D'$ ``sharing''
a point from $E$.} The following lemma is elementary:

\begin{lem}\label{trivial}
Let $D, D'$ be tiles in $A$ and  $E\subset A$.

1. If $I(D, D', E)\ge 4$, then $D=D'$.

2. If $I(D, D', E)\ge 3$,  then either $D\subset D'$ or $D'\subset
D$.
\end{lem}

In view of this lemma, we will say that a tile $D\in N$ that has a
thick vertex of $X$ (or an extreme point of $|N|$) on every edge
is {\em rigid} and a tile that has such a point on all but one
edges is {\em semi-rigid}. In particular, as $X$ is locally
finite, a compact subset of $X$ contains only finitely many rigid
tiles.

\begin{figure}[tbh]
\centerline{\epsfxsize=4.5in \epsfbox{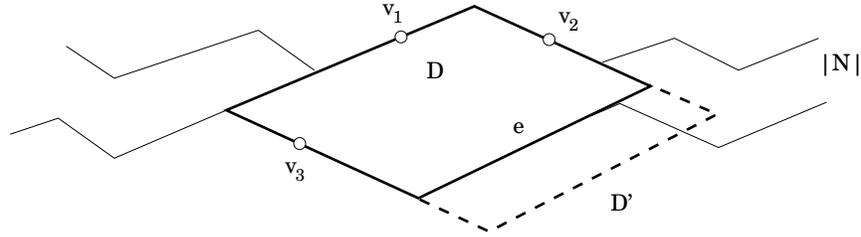}} \caption{\sl The
edge $e$ is not contained in $\D C$ and, hence, the tile $D$ can
be expanded to a tile $D'$. The points $v_1, v_2, v_3$ are the
thick vertices on $\D D \cap \D C$.} \label{f3}
\end{figure}

\begin{lem}\label{edge}
Let $D$ be a  tile in a weakly maximal corridor $N$, $|N|=C$, and
$e$ is an edge of $D$ which contains neither thick vertices, nor
extreme points of $C$. Then $e\subset \D C$.
\end{lem}
\proof If not, we can expand $D$ slightly to the direction of $e$
to a new tile $D'$ (which shares with $D$ three other sides), see
Figure \ref{f3}. Then $N\setminus \{D\} \cup \{D'\}$ is still a
tiled corridor in $X$ (which has the same extreme points as $N$)
contradicting weak maximality of $N$. \qed

\medskip
Suppose that $N$ is a weakly maximal tiled corridor in $X$,
$C=|N|$. Let $E$ denote the union of the set of extreme points
$\{v_-, v_+\}$ of $N$ with the set of thick vertices contained in
$B_d(\ol{v_- v_+})$, where $d=d(v_-, v_+)$. Let $t$ denote the
(finite)  cardinality of $E$. Set $E_C:=E\cap C$.

\begin{lem}
1. If $D, D'\in N$ are semi-rigid tiles so that $I(D, D', E_C)\ge
3$, then $D=D'$.

2. Suppose that $D, D'\in N$ are not semi-rigid tiles. If $I(D,
D', E_C)\ge 2$, then $D=D'$.
\end{lem}
\proof 1. Suppose $D\ne D'$. By Lemma \ref{trivial}, $D\subset D'$
or $D'\subset D$. Therefore, we can eliminate one of the tiles $D,
D'$ from the tiled corridor $N$ (without changing $|N|$),
contradicting the definition of a tiled corridor.

\begin{figure}[tbh]
\centerline{\epsfxsize=4.5in \epsfbox{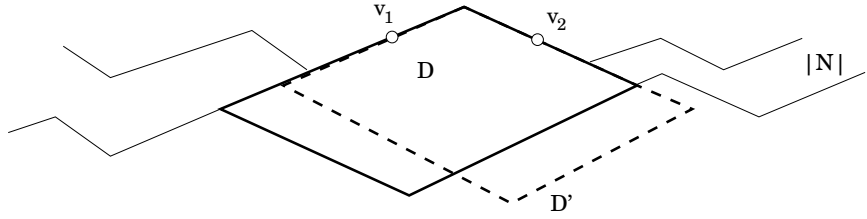}} \caption{\sl
$\D_-D$ is not contained in $\D C$. The points $v_1, v_2$ are the
thick vertices on $\D_+ D$.} \label{f2}
\end{figure}

2. Suppose that, say, $D, D'$ are not semi-rigid and share two
elements $a_1, a_2$ of $E$, both of which belong to $\D_+ C$.
Then, by Lemma \ref{edge}, $\D_- D, \D_- D'$ are both contained in
$\D_- C$. Suppose that $D\ne D'$. Definition of a tiled corridor
implies that $D$ cannot contain $D'$ and vice versa. But this
means that, $\D_- D$ has nonempty intersection with the interior
of $D'$, see Figure \ref{f2}. Contradiction.

The case when one of the shared points is on $\D_+ C$ and the
other on $\D_- C$ is similar and is left to the reader. \qed

The following three corollaries are then immediate:

\begin{cor}
Each tile in $N$ is uniquely determined by 2 or by 3 elements of
$E$.
\end{cor}

\begin{cor}\label{bound}
Let $N$ be a weakly maximal tiled corridor with the extreme points
$v_-, v_+$. Then number of tiles in $N$ is at most $t^3$.
\end{cor}

\begin{cor}\label{finiteness}
Given a tiled corridor $N$ in $X$, there exists only finitely many
weakly maximal tiled corridors $N'$ in $X$ so that $N\ll N'$.
\end{cor}

We then obtain:

\begin{prop}\label{maxexists}
For every tiled corridor $N_0$ there exists a maximal corridor $N$
so that $N_0\ll N$.
\end{prop}
\proof By Lemma \ref{wmax}, the set of weakly maximal tiled
corridors $N_i$ so that $N_0\ll N_i$ is nonempty. This set is
finite according to Corollary \ref{finiteness}. Therefore, it
contains a maximal element. \qed

\section{Combinatorial convexity}\label{convexity}

We continue with the notation of the previous section.  Our next
objective is to build convex subcomplexes in $X$ containing the
given pair of points:

\begin{thm}\label{convex}
Let $X$ be locally finite and modeled on a special Coxeter complex
$(A, W_{af})$, whose group of translations is countable and dense
in $\R^2$. Then:

1. For any two points $x_1, x_2\in X$, there exists a corridor $C$
in $X$ so that $x_1, x_2\in C$. Moreover, if $x_i$ is not a thick
vertex and $\ol{xy}$ is not a sub-wall near $x_i$, then we can
choose $C$ so that $x_i\notin \D C$.

2. For every corridor $C$ in $X$, there exists a convex corridor
$C'$ so that $C\subset C'$.
\end{thm}

The proof of this theorem will occupy the rest of the section.

\begin{rem}
1. Since $X$ contains only countably many walls, it contains only
countably many (convex) subcomplexes.

2. The above theorem is clear in the case when $X$ is a building,
as any two points belong to an apartment. In this case, we can use
$C=C'$, a single diamond.

3. If $X$ is not a building, the convex corridor $C'$ is, in
general, not a diamond. Accordingly, the weakly isometric
embedding $f: C'\to A$ will not send $C'$ to a convex subset of
$A$. See Remark \ref{nonconvex}.

4. The theorem is probably also true for arbitrary CAT(0) locally
Euclidean complexes (not necessarily 2-dimensional ones) provided
that we allow cells in $C$ to be arbitrary convex Euclidean
polyhedra. However, requiring $C$ in Theorem \ref{convex} to be
planar, imposes serious restrictions on the allowed subcomplexes.

5. In the main application of the above theorem, $X$ will be
finite, i.e., covered by finitely many apartments.
\end{rem}

We prove this theorem by:

i. Constructing a corridor $C$ containing $x_1$ and $x_2$. We will
subdivide the geodesic $\ga=\ol{x_1 x_2}$ into subsegments $\ga_i$
not passing through thick vertices. We then construct for each
$\ga_i$ a simple (possibly degenerate) corridor $C_i$, so that
$\ga_i\subset C_i$ and the end-points of $\ga_i$ are the extreme
vertices of $C_i$. The corridor $C$ will be the union of the
sub-corridors $C_i$. The most difficult part of the proof is
constructing a {\em planar neighborhood} $U_i$ of $\ga_i$. The
corridor $C_i$ will be contained in this neighborhood.

ii. Given a tiled corridor $N$ with $|N|=C$, we will take a
maximal tiled corridor $N'$ so that $N\ll N'$. In particular,
$C\subset C'=|N'|$. We then verify convexity of $C'$.

\medskip
{\bf 1. Constructing simple corridors.} We first assume that
$\ga=\ol{x_1 x_2}$ is a geodesic with distinct end-points, which
contains no thick vertices in its interior and is not contained in
a wall.

\begin{rem}
If $\ga$ is contained in a wall $L\subset X$, then $L\subset A$ is a planar set and
we obtain an isometric embedding
$$
\ga \embed L \embed A.
$$
Any sub-wall $e\subset L$ containing $\ga$, is therefore a
degenerate convex corridor in $X$ containing $\ga$.
\end{rem}

{\bf Planar neighborhoods.} At first, we assume, in addition, that
neither $x_1$, nor $x_2$ is a thick vertex and that $\ga$ does not
cross a thick wall. Then a {\em planar neighborhood} of
$\ga\subset X$ is a connected closed planar subset $\ol{U}$ of $X$
whose interior $U$ (in $X$) contains $\ga$.

\begin{lem}\label{openplanar}
Under the above assumptions, $\ga$ has a convex planar
neighborhood $\ol{U}$. In particular, an embedding $U\embed A$ will be distance-preserving.
\end{lem}
\proof Let $\eps$ be the minimal distance from $\ga$ to the union
of thick sub-walls of $X$. Our assumptions imply that $\eps>0$. We
then let $U$ be the open $\eps$-neighborhood $B_\eps(\ga)$ and
$\ol{U}$ be its closure in $X$. Then both $U$ and $\ol{U}$ are
convex in $X$, since $X$ is a CAT(0) space. Moreover, by the
choice of $\eps$,
$$
\ol{B_\eps(\ga)}=\bigcup_{x\in \ga} \ol{B_\eps(x)}
$$
where each $\ol{B_\eps(x)}$ is planar. Therefore, by considering the
developing map of $B_\eps$, and using convexity of this set we
obtain a morphism $f: B_\eps(\ga)\to A$. Convexity of
$U$ implies that $f$ is an isometric embedding. Therefore, it extends to an isometric embedding
$\ol{U}\to A$. \qed

\medskip
{\bf Relative planar neighborhoods of $\ga$.} We will have to
generalize the concept of a planar neighborhood to take into
account the cases when one of the $x_i$'s is a thick vertex, or
when $\ga$ crosses a thick sub-wall: In  both cases a planar
neighborhood of $\ga$ does not exist, so we will have to modify
the definition.

\begin{defn}
A {\em relative} planar neighborhood of $\ga$ is a closed connected
subset $\ol{U}\subset X$ (the closure of some $U\subset X$) so that the following holds:

1. $x_1, x_2\in \ol{U}, \ga \setminus \{x_1, x_2\}\subset U$.

2. If $x_i$ is not a thick vertex, then $U$ contains $x_i$.

3. If $x_i$ is a thick vertex, then (for sufficiently small
$\eps>0$) $S_i=\ol{U}\cap \ol{B_{\eps}(x_i)}$ is conical with
respect to $x_i$ and $\Si_x(S_i)$ is a chamber in $\Si_{x_i}(X)$.

4. There exists an injective isometric morphism $f: U\to A$ onto
an open subset of $A$. (Note that we do not require $U$ itself to
be open in $X$.)
\end{defn}

\begin{figure}[tbh]
\centerline{\epsfxsize=5in \epsfbox{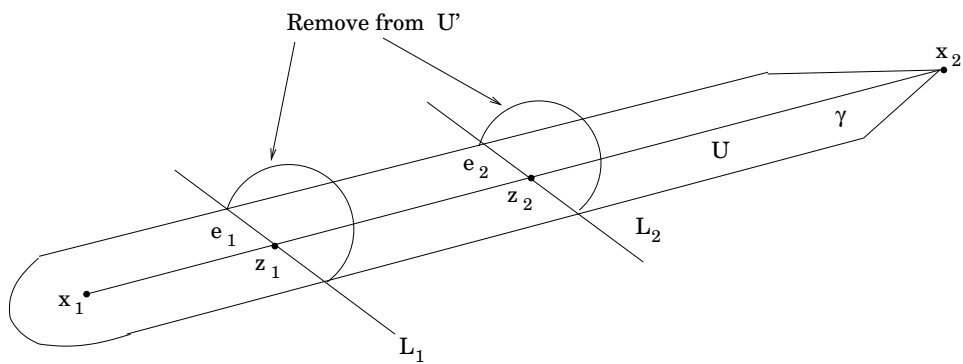}}
\caption{\sl Constructing a planar neighborhood.}
\label{f7}
\end{figure}

Figure \ref{f7} gives an example of relative neighborhood. We
observe that if $\ga$ is a geodesic satisfying the assumptions of
Lemma \ref{openplanar}, then the notions of a planar neighborhood
and of a neighborhood coincide.

\begin{lem}
Every geodesic $\ga=\ol{x_1 x_2}\subset X$ (which has no thick
vertices in its interior) has a convex  relative  planar
neighborhood.
\end{lem}
\proof First of all, without loss of generality we may assume that
for both $i=1,2$, either $x_i$ is a thick vertex or it is not
contained in a wall (otherwise, we can extend $\ga$ slightly in
the direction of $x_i$).

Let $\eps>0$ be such that $B_\eps(\ga)$ intersects only those
thick sub-walls which cross $\ga$, possibly at its endpoint(s). As
before, $B_\eps(\ga)$  is convex, however, it is no longer planar
(unless we are in the situation of Lemma \ref{openplanar}). Our
goal is to decrease this neighborhood of $\ga$ to make it planar
and to satisfy the rest of the properties of a relative
neighborhood.

Suppose that $x_i$ is a thick vertex. Consider the germ of $\ga$
at $x_i$.  Since $\ga$ is not contained in a wall, it follows that
this germ determines a regular point $\ga'(x_i)$ of
$\Si_{x_i}(X)$. Therefore, $\ga'(x_i)$ is contained in a unique
chamber $\si_i\subset  \Si_{x_i}(X)$. Hence, there exists a closed
conical (with respect to $x_i$) subset
$$
S_i\subset \ol{B_\eps(x_i)},
$$
so that
$$
\Si_{x_i}(S_i)=\si_i.
$$
Thus $S_i$ is isomorphic to an appropriate subset of
$$
Cone(\si_i).
$$
In particular, $S_i$ is planar. We then
remove from  $\ol{B_\eps(\ga)}$ the subset
$$
\ol{B_\eps(x_i)} \setminus S_i.
$$
If $x_i$ is not a thick vertex, we do no not modify $B_\eps(\ga)$
near $x_i$. This ensures that the resulting convex subset
$U'\subset B_\eps(\ga)$ satisfies Parts 1, 2, 3 of the definition.
We also observe that $U'$ were planar if the interior of $\ga$ did
not cross any thick sub-walls.

We now deal with the thick sub-walls which $\ga$ might cross. Let
$L_1,...,L_k$ denote the thick sub-walls of $X$ crossed by the
interior of $\ga$ at the points $z_1,...,z_k$. Then $U'$ is
locally planar except at the points of the geodesics $e_i=L_i\cap
U'$: Near those points, $U'$ is not even a  manifold, but looks
like an ``open book'' with the binding $e_i$. Our goal is to
``tear off'' all but two pages of this open book (for each
$i=1,...,k$), i.e., the ``pages''  which are not tangent to $\ga$.

Consider the link $\Si_{z_i}(X)$. Since $z_i$ is not a thick
vertex, this link is a building: This graph has exactly 2 thick
vertices, which are antipodal points representing the tangent
directions to $e_i$ at $z_i$. The tangent directions to $\ga$ at
$z_i$ are represented by two {\em regular} antipodal points in
$\Si_{z_i}(X)$. Therefore, there exists a unique apartment
$\si\subset \Si_{z_i}(X)$ containing these regular antipodal
points. We now remove from $U'$ the union of those open geodesic
segments $z_i u$ which are {\em not} tangent to $\si$.  (The
closure of the removed set is conical with respect to $z_i$.)  We
do this for each $z_i, i=1,...,k$. See Figure \ref{f7}. The
resulting set $U$ is still convex in $X$. Now, however, it is a
topological 2-manifold. We construct the injective morphism $f:
U\to A$ as in the proof of Lemma \ref{openplanar}. The fact that
the image of $f$ is open follows from the invariance of domain
theorem. \qed

\begin{cor}\label{simplecorridor}
For any pair of points $x_1, x_2\in X$ so that the open geodesic
segment  $x_1x_2$ is not a sub-wall and contains no thick
vertices, there exists a simple corridor $C\subset X$ containing
$x_1, x_2$. If $x_i$ is a vertex then we can choose $C$  that has
$x_i$ as one of its two extremes. If $x_i$ is not a thick vertex,
we can also choose $C$ so that $x_i\notin \D C$.
\end{cor}
\proof Let $\ol{U}$ be a convex (relative) planar neighborhood of
the geodesic $\ga=\ol{x_1 x_2}$, given by the above lemma.

1. We begin by considering the case when both $x_i$ are thick
vertices; thus, the germ of  $\ol{U}$ near $x_i$ is a sector $S_i$
(with the angle $\pi/m$ at $x_i$). Since $\ol{U}$ is planar, it
suffices to consider the case when it is contained in $A$.
Moreover, we can assume that it is embedded in $A$ so that the
slopes of the boundary of $S_1$ have the form $\pm \tan(\pi/2m)$.
Then, the same will be true for the slopes of $\D S_2$. The
geodesic segment $\ga\subset A$ is the graph of a linear function
$u: [a_1, a_2]\to \R$, with
$$
x_i=(a_i, u(a_i)), \quad i=1, 2.
$$
Now, approximate $u$ (from above and from below) by upper and
lower roof functions $u_\pm$ using Lemma \ref{approx}. Note that
the graphs of the germs of $u_\pm $ near $a_i$ are the boundary
segments of $S_i$, $i=1, 2$. By taking a sufficiently fine
approximation, we obtain that the graphs of $u_\pm$ are contained
in $\ol{U}$. Therefore, the region bounded by these graphs is the
required simple corridor.

2. Suppose that neither $x_1$ nor $x_2$ is a thick vertex and,
hence, $x_1, x_2\in U$. If the slope of $\ga\subset A$ represents
a vertex in $(S^1, W)$, we can as well assume that this slope
equals $\pm \tan(\pi/2m)$. In this case, using convexity of $U$,
it is easy to construct a single (``nearly degenerate'') tile
$D\subset U$ containing $x_1, x_2$. We therefore assume that the
slope of $\ga$ is in the open interval $(-\tan(\pi/2m),
\tan(\pi/2m))$. Pick disjoint congruent tiles $D_i\subset U$
containing $x_i, i=1, 2$. Let $y_i$ denote the extreme vertices of
the smallest tile $D\subset A$ containing $D_1\cup D_2$. Then
$y_i\in D_i, i=1, 2$.

\begin{rem}
If $x_i$ is a vertex, we can, if necessary, choose $D_i$ so that $y_i=x_i$.
\end{rem}

We assume that $a_1<a_2$ are the projections of $y_1, y_2$ to the
horizontal axis. The convex hull
$$
H=Hull(D_1\cup D_2)\subset U\subset A
$$
is a hexagon, with four edges $e^{\pm}_i\subset \D_\pm D_i$ incident to $y_i$,
$i=1, 2$  and the remaining two parallel edges $d^\pm$, so that
the boundary of $H$ is the union of two arcs
$$
e_1^+ \cup d^+ \cup e_2^+, \quad e_1^- \cup d^- \cup e_2^-.
$$
These arcs are the graphs of  functions $h^\pm: [a_1, a_2]\to \R$.
Approximating $h^\pm$ by the respectively upper and lower roof
functions $u_\pm$ we obtain a simple model corridor $C\subset U$
containing $x_1, x_2$. Since $C\subset A$ is a simple model
corridor, as a subcomplex in $X$ is automatically satisfies Axiom
C2.

3. The construction in the case when one of the points $x_1, x_2$ is a
vertex and the other is not, is a combination of  the  two arguments in Case 1 and 2,
and is left to the reader.  \qed

\medskip
{\bf 2. Constructing general corridors.} We now drop the
assumption that the open geodesic $x_1 x_2$ does not contain any
thick vertices.

\begin{prop}\label{ecorridor}
For every pair of points $x_1, x_2\in X$, there is a corridor $C$
containing $x_1, x_2$. If $x_i$ is not a thick vertex and the germ
of the geodesic $\ga=\ol{x_1 x_2}$ at $x_i$ is a regular direction
in $\Si_{x_i}(X)$, we can assume that $x_i$ does not belong to $\D
C$.
\end{prop}
\proof We subdivide the geodesic $\ga$ into maximal subsegments
$\ga_1\cup...\cup \ga_k$, so that  the interior of each
$\ga_i=\ol{y_i y_{i+1}}$ does not contain any thick vertices.
According to Corollary \ref{simplecorridor}, there exist simple
(possibly degenerate) corridors $C_1,...,C_k$ so that each $C_i$
contains $\ol{y_i y_{i+1}}$.  Moreover, for $1<i<k$, the vertices
$y_i, y_{i+1}$ are the extremes of $C_i$.  If $C_1$ or $C_k$ is
non-degenerate and the corresponding point $x_i$ in this corridor
is not a thick vertex, by Corollary \ref{simplecorridor}, we can
assume that  $x_i$ does not belong to the boundary of this
corridor.

Let $C$ denote the union $C_1\cup ...\cup C_k$. We claim that $C$
is the required corridor.

For each $i$ we have an isomorphic (and isometric) embedding $f_i:
C_i\to A$, whose image is a simple model corridor $C_i'\subset A$.
We combine these maps to obtain a weakly isometric embedding $f:
C\to A$ as follows: We choose $f_1$ arbitrarily. We then normalize
each $f_i$ (by a translation in $W_{af}$), so that
$$
y_i':=f_i(y_{i})=f_{i+1}(y_i)
$$
and $y_i'$ is the extreme rightmost point of $C_i'$ and extreme
leftmost point of $C_{i+1}'$. Since each $C_i'$ is a model
corridor, this ensures that the germs of $C_i', C_{i+1}'$ at
$y_i'$ are antipodal in $\Si_{y_i'}(A)$, provided that both are
chambers. In any case,
$$
\del=\pi - \pi/m \le \angle( \Si_{y_i'}(C_i'),
\Si_{y_{i+1}'}(C_{i+1}'))\le \pi
$$
with equality to $\del$, except for the case when both $C_i', C_{i+1}'$ are degenerate and their
union is a geodesic segment in $A$. In the latter case the angle in $A$ equals $\pi$ and the angle in
$\Si_{y_i}(X)$ is $\ge \pi$, since $\ga$ is geodesic.

We thus obtain a morphism $f: C\to A$ whose restriction to each
$C_i$ is $f_i$. For each $i$, the convex hull of $C_i'$ is a tile
$D_i'\subset A$ with the extreme vertices $y_i', y_{i+1}'$. By the
construction, for $i\ne j$, $D_i'\cap D_j'$ is at most one point
(the common extreme vertex). Therefore, $f$ is injective.

We next need to verify the weak isometry at the vertices $y_i$.
Suppose that both $C_i, C_{i+1}$ are nondegenerate corridors. Then
their germs at $y_i$ are chambers $\si_i, \si_{i+1}\subset
\Si_{y_i}(X)$. Since $\ga$ is a geodesic, its angle at $y_i$ is
$\ge \pi$. Therefore, the chambers $\si_i, \si_{i+1}$ are {\em
super-antipodal}: The angular distance between these chambers in
$\Si_{y_i}(X)$ is $\ge \del=\pi- \frac{\pi}{m}$. On the other
hand, the germs of $C_i'$ at $y_i'=f(y_i)$ are antipodal in
$\Si_{y_i'}(A)=S^1$. Therefore, the map
$$
f: \Si_{y_i} C \subset \Si_{y_i} X \to \Si_{y_i'}(A)
$$
is a  $\del$-isometry. Thus, $f: C\to A$ is a weakly isometric embedding.

Suppose that $C_i$ is degenerate, but $C_{i+1}$ is not. The germ $\si_i$ represents the germ of
the geodesic $\ga_i$ at $y_i$ and $\si_{i+1}$ is a chamber containing the germ of the geodesic $\ga_{i+1}$
at $y_i$. Therefore, since $\si_{i+1}$ contains a point antipodal to $\si_i$, the angular distance in
$\Si_{y_i}(X)$ between the vertex $\si_i$ and the chamber $\si_{i+1}$ is $\ge \del$. The angular distance
between their images $f(\si_i), f(\si_{i+1})$ in $S^1$ is exactly $\del$. Therefore,
$f$ is a  weak isometry at $\Si_{y_i}(C)$.

The argument when both $C_i, C_{i+1}$ are degenerate follows from our discussion of the angles between
$C_i, C_{i+1}$ and $C_i', C_{i+1}'$.

Thus, $C$ satisfies Axiom C1 of corridors. We need to check Axiom
C2. Suppose that $C$ admits a bridge
$$
D_0\cup D_1\cup \cup D_2 \cup T,
$$
where $D_0=\ol{v_1v_2}$ is degenerate. Then, by the construction
of the corridor $C$, one of the $v_1, v_2$ (say, $v_1$) is a
vertex $y_i$ and $D_0\subset C_{i+1}$ (one can actually assume
that $D_0=C_{i+1}$). Therefore, the germ of the geodesic $\ga$ at
$v_1$ is contained in $D_1\cup D_0$. Since $T\cup D_1$ is a planar
region, this means that the angle of $\ga$ at $v_1=y_i$ is less
than $\pi$, contradicting the fact that $\ga$ is a geodesic. \qed

\begin{lem}\label{maxconvex}
If $N$ is a maximal tiled corridor (with respect to the order
$\ll$) in $X$, then $C=|N|$ is convex in $X$.
\end{lem}
\proof Suppose that  $C$ is not convex in $X$. Then, since $X$ is
CAT(0), there exists a vertex $v$ of $C$ so that $C$ fails to be
locally convex at $v$. Since all the angles of $C$ are of the form
$\pi/m, \pi+\pi/m$, it follows that $C$ has the angle $\pi+\pi/m$
at $v$. This means that $\Si_v(C)$ is not convex in $\Si_v(X)$. In
other words, if $\eta, \xi$ denote the tangent directions to $\D
C$ at $v$, then $\angle(\eta,\xi)< \pi + \pi/m$. Recall that
$\Si_v(X)$ is a bipartite graph with the edges of the length
$\pi/m$. Since $\xi_1, \xi_2$ have different colors, it follows
that $\angle(\xi_1,\xi_2)\le \pi- \pi/m$. Since $\Si_v(X)$ has
girth $\ge m$, it follows that $\angle(\xi_1,\xi_2)= \pi- \pi/m$.
Therefore, there exists a small flat sector $S\subset X$ with the
vertex $v$ and $\Si_v(S)=\ol{\xi_1\xi_2}\subset \Si_v(X)$. Thus,
by density of the walls of $A$, we can find a small tile $D\subset
S$  with the vertex $v$ and the angle $\pi- \pi/m$ at $v$. The
union $C'=C\cup D$ is, therefore, again a corridor and $N'=N \cup
\{D\}$ is a tiled corridor so that $N\ll N'$. Since $|N'|$ is
strictly larger than $C$, it follows that $N$ is not maximal.
Contradiction. \qed

\begin{rem}
In general, weakly maximal corridors (i.e., maximal with respect
to the order $<$) need not be convex.
\end{rem}

\begin{cor}\label{convexexist}
Every corridor $C\subset X$ is contained in a convex corridor
$C'\subset X$.
\end{cor}
\proof Let $C=|N|$. Let $N'$ be a maximal corridor so that $N\ll
N'$ (see Proposition \ref{maxexists}). Then $C\subset C'=|N'|$.
The corridor $C'$ is convex by the above lemma. \qed

\begin{rem}
Although we do not need this, given $N$ so that $C=|N|$, we can
find a canonical convex corridor $C'$ containing $C$. Namely, the
set of weakly maximal tiled corridors $N_i$ so that $N\ll N_i$ is
finite according to Corollary \ref{finiteness}. The number of
maximal corridors among the corridors $N_i$ is finite as well. Let
$C'$ denote the intersection
$$
C'=\bigcap_i |N_i|
$$
taken over all maximal corridors $N_i$. One can verify that $C'$
is indeed a corridor. Convexity of $C'$ follows from the fact that
it is intersection of convex sets.
\end{rem}

\medskip
We can now finish the proof of Theorem \ref{convex}.

\proof  The first assertion of Theorem follows from Proposition
\ref{ecorridor}. The second assertion is  Corollary
\ref{convexexist}. \qed

\section{Generalizations of the convexity theorem}\label{generalized}

In this section we prove two generalizations of Theorem
\ref{convex}. Our first goal is to provide a generalization which
allows corridors to have infinite diameter.

We first have to generalize the notion of a {\em tile}. We define
{\em extended (infinite) tiles}, which are infinite sectors
bounded by walls with the slopes $\pm \tan(\pi/2m)$ and having the
vertex angle $\pi/m$ or $\pi-\pi/m$. We will also allow the whole
apartment as a single tile. Such infinite tiles appear as limits
(in the Chabauty topology on the space of closed subsets of $A$)
of sequences of the ordinary tiles. Accordingly, we generalize the
upper and lower roof functions to allow functions defined on $\R$.
We still require such piecewise-linear functions to have only
finitely many break-points: Every roof function admits a partition
of $\R$ in finitely many subintervals, where the function is
linear. Then, a {\em extended model corridor} is the closed region
bounded from above and below by graphs of upper and lower roof
functions; we also declare $A$ to be a corridor (consisting of a
single infinite tile). It is clear that each corridor $C\subset A$
is a union of finitely many tiles, two of which are infinite, with
the vertex angles $\pi/m$. Therefore, for every corridor $C$,
there exists a extended tiled corridor $N$ so that $|N|=C$: The
corridor $N$ still has only finitely many tiles, but some of them
are allowed to be infinite.

Thus the ``extremes'' of the extended corridors $C$ are points at
infinity of $C$. Namely, if $\ga\subset A$ is any complete
geodesic contained in $C$, we declare the points at infinity of
$\ga$ to be the extremes of $C$. Therefore, $C$ now has
(typically) more than 2 extreme points.

Let $X\in \c$, a space modeled on $(A, W_{af})$ and covered by
finitely many apartments. (It is no longer enough to assume local
finiteness.) We define (infinite) corridors in $X$ by repeating
the definition for the ordinary corridors.

\begin{prop}
Every complete geodesic $\ga\subset X$ is contained in a corridor.
\end{prop}
\proof The arguments remain pretty much the same as in the proof
of Proposition \ref{ecorridor}, except that we start by observing
that, for sufficiently large $t$, the subrays $\ga(-\infty, -t),
\ga(t, \infty)$ in $\ga$ do not cross any walls in $X$. Therefore,
there are infinite tiles $D_1, D_2\subset X$ containing these
subrays.

\begin{figure}[tbh]
\centerline{\epsfxsize=5in \epsfbox{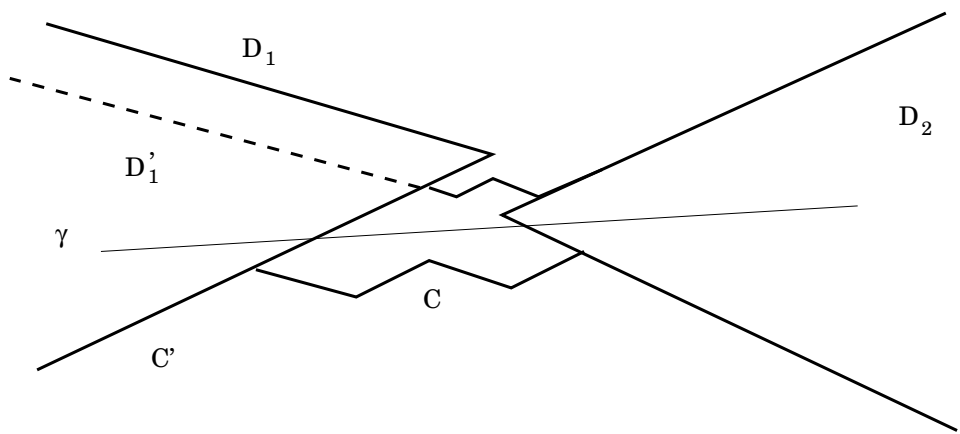}} \caption{\sl
Constructing an extended corridor.} \label{f6}
\end{figure}

Then $\ga \setminus (D_1 \cup D_2)$ is a finite geodesic segment
and we can cover it with a tiled corridor $C$ using Proposition
\ref{ecorridor}. One then shrinks the tiles $D_1, D_2$, to $D_1',
D_2'$ so that
$$
C'=D_1'\cup C\cup D_2'
$$
is a corridor covering $\ga$, see Figure \ref{f6}. \qed

\medskip
The definitions of the orders $<$ and $\ll$ on corridors
generalize almost verbatim to the corridors having infinite tiles.
The only modification one needs is that instead of dealing with
finite extreme vertices, we require an extended corridors to
contain certain points at infinity. Equivalently, one can ask for
corridors covering a fixed complete geodesic in $X$.

One again sees that for every extended tiled corridor $N$, there
exists an extended weakly maximal tiled corridor $N'$ so that
$N<N'$: While considering an infinite increasing sequence of tiled
corridors we no longer can argue that the tiles have uniformly
bounded diameter, but such unbounded sequences of tiles still
converge to infinite tiles. (The limit could be the whole
apartment $A$.) Given a weakly maximal extended corridor $N$ we
construct a maximal extended corridor $N'$ so that $N\ll N'$ by
appealing to the fact that  $X$ has only finitely many thick
vertices.

The proof that every maximal extended corridor is convex is
exactly the same as that one of Lemma \ref{maxconvex}. By
combining these observations, we obtain

\begin{thm}\label{convex1}
1. Every complete geodesic $\ga \subset X$ is contained in an
extended corridor $C$.

2. Every extended corridor $C$ is contained in a maximal extended
corridor $C'$. The corridor $C'\subset X$ is convex.

3. Let $\De_i, i=1, 2$ be Weyl sectors which share the tip $x\in
X$ and which are super-antipodal at $x$. Then there exists an
extended corridor $C\subset X$ so that
$$
\De_1\cup \De_2\subset C.
$$
Moreover, the corridor $C$ contains a geodesic with the ideal
points in $\geo \De_i, i=1,2$.
\end{thm}
\proof It remains to prove Part 3. We start with the case when $x$
is a vertex and, hence, $\De_i$ are bounded by sub-walls. Then,
$\De_i, i=1, 2$, are extended tiles  in $X$ and, by our
assumptions, $C=\De_1\cup \De_2$ is an extended corridor. The
geodesic contained in this corridor is the union of rays $\rho_i$
which are super-antipodal at $x$ and asymptotic to points in
$\De_i, i=1, 2$.

If $x$ is not a vertex, $\De_1, \De_2$ have to be antipodal at
$x$. Moreover, the boundary rays of $\De_{i}$ are not sub-walls.
In particular, they do not contain thick vertices. Therefore,
there exist (special) vertices $x_i\in \De_{i+1}$ close to $x$ and
Weyl sectors $\De'_i$ with the tips $x_i$, so that
$$
\De_i\subset \De_i', i=1, 2
$$
and
$$
\De_1'\cap \De_2'=D'
$$
is a tile. Therefore, we obtain an extended corridor $C=\De_1'\cup
\De_2'$ containing $\De_1\cup \De_2$. \qed

\medskip
We now return to the discussion of ordinary corridors. In the
following theorem it suffices to assume that $X$ is locally
finite, modeled on $(A, W_{af})$.

\begin{thm}\label{convex2}
For $x_1, x_2\in X$ let $\si_i$ be chambers in the buildings
$\Si_{x_i}(X)$, which are given structure of thick buildings,
$i=1,2$. Then there exists a convex corridor $C\subset X$
containing $x_1, x_2$, so that
$$
\si_i\subset \Si_{x_i}(C), i=1, 2.
$$
\end{thm}
\proof As before, set $\ga=\ol{x_1 x_2}$. We start with a convex
corridor $C$ containing $x_1, x_2$, given by Theorem \ref{convex}.
We modify $C$ twice: First to ensure
$$
\si_1\subset \Si_{x_1}(C)
$$
and next, to ensure
$$
\si_2\subset \Si_{x_2}(C).
$$
We describe only the first modification, since the second is
obtained by applying the same procedure to the modified convex
corridor and switching the roles of $x_1, x_2$. To simplify the
notation, set $x:=x_1$, $\si:=\si_1$. We let $\eta$ be the tangent
direction to $\ga$ at $x$.

\medskip
Case 1 (the generic case): $x$ is not a thick vertex, is not
contained in a thick sub-wall and the germ of $\ga$ at $x$ is not
contained in a sub-wall. In particular, $\si= \Si_x(X)$ is the
apartment. Then, according to Theorem \ref{convex}, we can assume
that $x\notin \D C$. Therefore $\si=\Si_x(X)=\Si_x(C)$.

Case 2: We now allow $x$ to be contained in a thick sub-wall, but
retain the rest of the restrictions. Since $x$ is not a thick
vertex, $\Si_x(X)$ is a building (but not a single apartment). The
chamber $\si$ is a half-apartment in $\Si_x(X)$ and $\eta$ is a
regular direction.

(a) If $\eta\in \si$, then $\si=\Si_x(C)$ and we are done.

(b) If $\eta\notin \si$, then, there exists a (necessarily
regular) direction $\xi\in \si$ antipodal to $\eta$ and we can
replace $\ga$ with a slightly longer geodesic $\t\ga=\ol{\t{x}
{x_2}}$ tangent to $\xi$. The convex corridor $\t{C}$ containing
$\t{\ga}$  will then contain $x$ in its interior and, therefore
$$
\ga\subset \t{C}, \si\subset  \Si_x(\t{C}).
$$
Hence, we are done in this case as well.

Case 3:  $x$ is a thick vertex, but $\eta$ is a regular direction,
i.e., the germ of $\ga$ at $x$ is not contained in a sub-wall.

(a) Suppose that
$$
\angle(\xi, \eta)\le \pi, \quad \forall \xi\in \si.
$$
This means that $\si$ is contained in an arc $\tau=\ol{\eta
\zeta}\subset \Si_x X$ of the length $\le \pi$. Since $\Si_x(X)$
has extendible geodesics, we can assume that the length of $\tau$
is $\pi$ and, hence, $\zeta$ is antipodal to $\eta$. We now use
the same arguments as in case 2a: Extend $\ga$ to a longer
geodesic $\t{\ga}$ tangent to $\zeta$. Let  $\t{C}$ be the convex
corridor containing $\t\ga$. By convexity, $\Si_x(\t{C})$ will
contain a geodesic in $\Si_x(X)$ connecting $\eta, \zeta$. Since
$\eta$ is regular, this geodesic is unique and equals $\tau$. In
particular, $\si\subset \Si_x(\t{C})$.

(b) If (a) fails then $\si$ contains a regular direction $\xi$ so
that
$$
\angle(\xi, \eta)\ge \pi.
$$
We again repeat the arguments in the case 2(b): Take a longer
geodesic $\t{\ga}$ tangent to $\xi$. Then the corresponding convex
corridor $\t{C}$ will contain a germ of $\si$ at $x$.

Case 4: $\eta$ is a singular direction and the germ of $\ga$ at
$x$ is contained in a sub-wall. Then $C$ contains a degenerate
irreducible component $C_1$ so that $x\in C_1$. We have the same
dichotomy as in Case 3:

\begin{figure}[tbh]
\centerline{\epsfxsize=5in \epsfbox{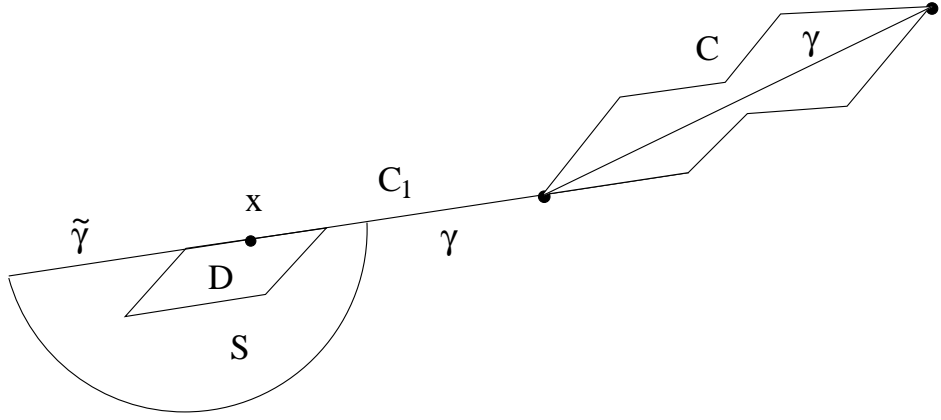}} \caption{\sl
Constructing the corridor $\t{C}=D\cup C$.} \label{f10}
\end{figure}

(a) Suppose that $\si$ is contained in an arc $\tau=\ol{\eta
\zeta}\subset \Si_x X$ of the length $\pi$. We construct $\t\ga$
as before. There exists a flat sector $S$ (germ of a
half-apartment) with the tip $x$ so that $\tau= \Si_x(S)$. We then
find a tile $D\subset S$ so that $\tau=\Si_x(D)$, i.e., one of the
edges of $D$ is contained in $\t\ga$, so that $C_0=D\cup C_1$ is a
corridor. See Figure \ref{f10}. Set
$$
\t{C}:=D \cup C.
$$
By the construction, $\si\subset \tau= \Si_x(D)= \Si_x(\t{C})$. We
claim that $\t{C}$ is a corridor. We obtain a weakly isometric
embedding $\t{C}\to A$ by combining the weakly isometric
embeddings of $C_0$ and $C$. If $\t{C}$ were to admit a bridge,
there would have to be at least two vertices where $\t{C}$ fails
to be locally convex. One of them is in $C_0$. Since $C$ is
convex, it does not have any nonconvex vertices. Thus, $\t{C}$ is
a corridor.

(b) If $\si$ contains a regular direction $\xi$ so that
$\angle(\xi \eta)\ge \pi$, we repeat the arguments in Case  3(b).
\qed

\begin{cor}
Let $P_j, j=1, 2$ be (convex compact) cells in $X$. Then there
exists a finite set of corridors $C_1,...,C_k$ in $X$ so that for
every $(x_1, x_2)\in P_1\times P_2$ and chambers
$$
\si_j\subset \Si_{x_i}(X), j=1, 2,
$$
there exists a corridor $C_i, i=1,...,k$ containing $x_1, x_2$ and
so that
$$
\si_j\subset \Si_{x_i}(X), i=1, 2.
$$
\end{cor}
\proof The product $\Si(P_1)\times \Si(P_2)$ is compact. For every
subcomplex $C$ (we are interested only in corridors, of course),
$$
\Si( int_{P_i}(C\cap P_i))\subset \Si( P_i)
$$
is an open subset, since $P_i$ is contained in an apartment in $X$.
Therefore, we obtain an open
covering of  $\Si(P_1)\times \Si(P_2)$ by the sets of the form
$$
\Si( int_{P_i}(C\cap P_i)), \quad C \quad \hbox{ is a convex corridor in~~~ } X 
$$
Therefore, by compactness, there exists a finite number of
corridors $C_1,...,C_k$ so that for every $x_i\in P_i$ and every
chamber $\si_i\in \Si_{x_i}(P_i)$,
$$
x_i\in C, \quad \si_i\subset \Si_{x_i}(C). \qed
$$

\section{Enlargement along convex corridors}\label{category}

For a space $X\in \c$ we define $Corr(X)$, the set of all {\em
convex extended corridors} in $X$. In this section we will
frequently omit the adjective {\em extended} since almost all the
discussion of corridors will be local.

\begin{defn}
Given $X\in \c$, we define an {\em enlargement of $X$ along a
(extended) convex corridor $C\subset X$} to be the space $X'=
X_C=X\cup_C A\in \c$ obtained by attaching to $X$ the apartment
$A$ along $C\subset X$. The attachment is given by a weakly
isometric embedding $f: C\embed A$. We will also use the term {\em
enlargement} for the tautological embedding $\iota: X\to X'$.
\end{defn}

It is important to note that every wall in $X$ is an extended
corridor, therefore, enlargements include attaching apartments
along walls.

\begin{lem}
If $X'$ is obtained from $X$ by enlargement, then $X'$ is CAT(0).
\end{lem}
\proof Since $C$ is contractible, $X'$ is simply-connected. Since
$X'$ is a locally-Euclidean cell complex, it suffices to check
that the links of vertices $v$ of $X$ are CAT(1) spaces. Let $v$
be a vertex of $C$. If $f(C)$ is locally convex at $f(v)$ then
$$
\Si_v(X')=\Si_v(X)\cup_Z \Si_{f(v)}(A)
$$
is again CAT(1), since $Z=\Si_v(C)$ is convex in both
$Y=\Si_v(X)$ and $\Si_{f(v)}(A)$.

Consider therefore, a vertex $f(v)$ where $f(C)\subset A$ is not
locally convex. Then $Y'=\Si_v(X')$ is obtained from $Y=\Si_v(X)$
by attaching either one arc $\al$ (if $v$ is not a cut-vertex of
$C'$) or two arcs $\al_1, \al_2$ (if $v$ is a cut-vertex); in
either case, the length of $\al$ (and $\al_1, \al_2$) is greater
than or equal to $\del=\pi(1-\frac{1}{m})$. We will consider the
case of one arc since the case of two arcs is similar. Since
$C\subset X$ is convex, the arc $\al$ is attached to vertices
$\xi_1, \xi_2$ of $Y$, so that $\angle_Y(\xi_1,\xi_2)\ge
\pi(1+\frac{1}{m})$ (which is the inner angle of $C$ at $v$). Let
$\si\subset Y'$ be the shortest embedded circle. If $\si$ does not
contain $\al$, then $\si\subset Y$ and $length(\si)\ge 2\pi$. If
$\si$ contains $\al$, then $\si=\al\cup \be$, where $\be$ is an
arc in $Y$ connecting $\xi_1, \xi_2$. Then, $length(\be)\ge
\pi(1+\frac{1}{m})$ and, hence
$$
length(\al)\ge  \pi(1-\frac{1}{m}) + \pi(1+\frac{1}{m})=2\pi.
$$
Thus, $Y'$ is CAT(1) and $X'$ is CAT(0). \qed

\medskip
We therefore obtain a category $\e$, whose objects are element of
$\c$ (i.e., finite spaces modeled on $(A,W_{af})$) and whose
morphisms, $Mor(\e)$,  are compositions of enlargements.

\begin{lem}\label{conserve}
Every enlargement $\iota: X\to X'=X_C$ is a weakly isometric
embedding. The map $\iota$ is 1-Lipschitz, in particular, its
restriction to every apartment in $X$ is an isometry.
\end{lem}
\proof Verifying that a map is a weakly isometric embedding is a local problem.
Let $v\in X$ be a vertex. Then $Y'=\Si_v(X')$ is obtained from $Y=\Si_v(X)$ by attaching at
most two arcs, each having length $\ge \del=\pi  - \frac{\pi}{m}$.
Therefore, the inclusion $Y\embed Y'$ is a $\del$-isometric
embedding. It follows that the map $X\embed X'$ is a weak
isometry.

The fact that $\iota$ does not increase distances is clear. Suppose that
$\si\subset \Si_x(X)=Y$ is an embedded  circle of length $2\pi$
(say, an apartment). Set $Y':= \Si_x(X')$.
Suppose that there exist $\xi, \eta\in \si$ so that
$$
\angle_{Y} (\xi, \eta) > \angle_{Y'} (\xi, \eta).
$$
Then, combining the geodesic $\ol{\xi \eta}\subset Y'$ with the
geodesic $\ol{\xi \eta}\subset Y$, we obtain an embedded circle of
length $<2\pi$ in $Y$. This contradicts the fact that $X'$ is
CAT(0) and $Y'$ is CAT(1). Therefore, for every apartment
$A'\subset X$, its image $\iota(A')$ is convex in $X'$ and, hence,
the restriction of $\iota$ to $A'$ is an isometry.  \qed

\medskip
Note that $\iota: X\to X'$ is very seldom an isometric embedding,
namely, if and only if $f(C)\subset A$ is convex, which means that
$C$ is a single (extended) tile.

\begin{lem}\label{corridorconservation}
Let $X'=X\cup_{C'} A$ be obtained from $X$ by enlargement along a
(extended) convex corridor $C'\subset X$ and $\iota: X\to X'$ be
the tautological embedding. Let $C\in Corr(X)$. Then the image
$\iota(C)\subset X'$ is again a corridor (not necessarily convex).
\end{lem}
\proof  Proof of Axiom C1 is the same as in Lemma \ref{isomembed}.
We need to verify Axiom C2. We will identify $C$ and $\iota(C)$.
Suppose that $C\subset X'$ does admit a bridge $D_1\cup D_2 \cup D_0 \cup T$. Let $e_i\subset \D D_i$
be the edges adjacent to the trapezoid.

Since $C$ admits a bridge, it is not convex in $X'$. Since $C$ was
convex in $X$, this means that $\iota$ is not a local isometry at
the vertices $v_1, v_2$ of the bridge (where $D_0=\ol{v_1 v_2}$).
It follows that $v_1, v_2$ are vertices of $\D C'$. By convexity,
$D_0\subset C'$ as well. Moreover, the germs of $e_i \cup D_0$ at
$v_i$ have to be contained in $\D C'$ ($i=1, 2$), otherwise,  they
would be convex in $X'$. Thus, $D_0\subset \D C'$ and, therefore,
$C'\subset X'$ admits a bridge with trapezoid $T$. Since the
convex hull of $C'\subset X'$ is contained in the apartment $A$
(which was attached to $X'$ along $C'$), we conclude that
$T\subset A$. This, however, contradicts the fact that the image
of $C'$ in $A$ is a model corridor. \qed

\begin{cor}\label{largercorridor}
Let $\iota: X\to X'$ be a morphism in $\e$. Then for every $C\in
Corr(X)$, there exists a (extended) convex corridor $C'\subset X'$
containing $\iota(C)$.
\end{cor}
\proof Since $\iota$ is a composition of enlargements, we
construct $C'$ inductively starting with $C_1=C$. If $C_i\subset
X_i$ is a convex corridor and $\iota_i: X_i \to X_{i+1}$ is an
enlargement, then $\iota_i(C_i)\subset X_{i+1}$ is a corridor
according to Lemma \ref{corridorconservation}. We then replace the
corridor $\iota_i(C_i)$ with a convex corridor in $X_{i+1}$
containing $\iota_i(C_i)$ (see Part 1 of Theorem \ref{convex1}).
\qed

\section{Construction of the building}\label{mainproof}

We are now ready to prove Theorem \ref{main}. We proceed by
constructing inductively an increasing sequence $X_n\in \c$, $n\in
\N$, by using enlargements along (extended) convex
corridors\footnote{Recall that every wall in $X_n$ is in
$Corr(X_n)$.} $C\in Corr(X_n)$.  Then the building $X$ will be the
direct limit of the sequence $(X_n)$.

Set
$$
Corr:= \bigsqcup_{n\in \N} Corr(X_n),
$$
Since $Corr(X_n)$ is countable for every $X_n\in \c$, we get a
bijection $\om_n: \N\to Corr(X_n)$. Define
$$
\om: \N^2\to Corr
$$
by
$$
\quad \om(m,n)=\om_n(m).
$$

Let $\nu=(\nu_1, \nu_2): \N\to \N^2$ denote the standard enumeration, so that both $\nu_i$ are
1-Lipschitz functions $\N\to \N$. In particular,
\begin{equation}\label{lip}
\nu_2(n)\le n, \quad \forall ~n\in \N.
\end{equation}
Lastly, set $\Om:=\om\circ \nu$.

We start with $X_1=A$, the model apartment. However, we could as
well start with any $X_1\in \c$.

Then, for  $n\in \N$ we compute the (extended) convex  corridor
$C\in Corr(X_k)$, $C=\Om(n)$, where $k=\nu_2(n)\le n$. We have a
morphism $\iota: X_k \to X_n$ which is a composition of
enlargements.   According to Corollary \ref{largercorridor}, there
exists a convex (extended) corridor $C'\in Corr(X_n)$ containing
$\iota(C)\subset X_n$. The corridor $C'$ is not unique. This does
not matter with one exception: If $C$ is a geodesic in an
apartment in $X_k$, $\iota(C)$ is still a geodesic in $X_n$; in
this case  we take $C':=C$.

We then  set
$$
X_{n+1}=X_{n,C'}=X_n \cup_{C'} A
$$
be the enlargement of $X_n$ along $C'$.

To make sure that $X_{n+1}$ is well-defined, we need that $\Om(n)$
is a corridor in $X_k$, for some $k\le n$. This immediately
follows from the inequality (\ref{lip}).

Our next task is to give $X$ structure of a building. Observe that
each apartment in $X_n$ is also an apartment in $X_{n+1}$ (see
Lemma \ref{conserve}). Therefore, we define the set of apartments
in $X$ as the union of the sets of apartments in $X_n$, $n\in \N$.
Since each $X_n$ is modeled on $(A, W_{af})$, we conclude that $X$
is also modeled on $(A, W_{af})$: The transition maps are in
$W_{af}$ and for every pair of apartments $\phi, \psi: A\to X$
$$
\phi^{-1} \psi(A)\subset A
$$
is a closed and convex subset in $A$ bounded by sub-walls.

\begin{lem}
Let $x_i\in X, i=1, 2$ and $\si_i\in \Si_{x_i}(X)$ be germs of Weyl chambers.
Then there exists an apartment
in $X$ containing $x_1, x_2$ and the germs $\si_1, \si_2$. In particular,
every two points of $X$ belong to a common apartment.
\end{lem}
\proof Take the smallest $k$ so that $x_i\in X_k$ and
$\si_i\subset \Si_{x_i}(X_k)$, $i=1,2$. According to Theorem
\ref{convex2}, there exists a convex corridor $C\subset X_k$
containing $x_i$ and $\si_i, i=1, 2$. Then there exists $n\ge k$
so that $C=\Om(n)$. Therefore, $X_{n+1}$ is obtained by attaching
an apartment $A$ along a convex corridor $C'\subset X_n$
containing $\iota(C)$, $\iota: X_k \to X_n$. It follows that the
images of $x_i, \si_i$ in $X_{n+1}$ will belong to an apartment
$A$. Hence, the image of $A$ in $X$ will contain $x_i, \si_i$,
$i=1,2$. \qed

\medskip
We can now give $X$ structure of a CAT(0) space, so that every
apartment is isometrically embedded.  Even though the morphisms
$\iota: X_k \to X_n$ are not isometric, nevertheless,  for every
$x_1, x_2\in X_k$ there exists $n$ so that $x_1, x_2$ belong to an
apartment in $X_n$; therefore, for every $l\ge n$ we have:
$$
d_{X_l}(x_1, x_2)= d_{X_n}(x_1, x_2).
$$
We, hence, obtain a metric on $X$. Since every morphism in $\e$ is
1-Lipschitz, it follows that for every $n$, $X_n \to X$ is also
1-Lipschitz. The metric on $X$ is a direct limit of the metrics on
the CAT(0) spaces $X_n$; therefore $X$ is also CAT(0) (as the
CAT(0) condition is a condition on quadruples of points).

\begin{lem}
Let $\De_1, \De_2\subset X$ be antipodal Weyl sectors with the
common tip $x$. Then there exists an apartment $A\subset X$
containing $\De_1\cup \De_2$.
\end{lem}
\proof We first find $k$ so that $\De_1\cup \De_2\subset X_k$.
Since $\De_1, \De_2$ are antipodal in $X$, and $X_k \to X$ is
1-Lipschitz, it follows that $\De_1, \De_2$ are super-antipodal in
$X_k$. Therefore, according to Theorem \ref{convex1}, there exists
an extended convex corridor $C\subset X_k$ containing $\De_1\cup
\De_2$. The corridor $C$ equals $\Om(n)$ for some $n\ge k$. Recall
that $X_{n+1}$ is obtained from $X_n$ by enlargement along a
convex corridor containing the image $\iota(C), \iota: X_k\to
X_n$. It follows that the images of $\De_1, \De_2$  in $X_{n+1}$
(and, hence, in $X$) are contained in an apartment.  \qed

We thus verified axioms A1$^\prime$ and A2$^\prime$ of Euclidean
buildings for $X$.

\begin{lem}
$X$ is thick.
\end{lem}
\proof Let $L$ be a wall in $X$. Then $L$ is contained in an
apartment $A\subset X$, which means that $L$ is a wall in an
apartment $A\subset X_k$ for some $k$. Since $L$ is also a
(degenerate) corridor it follows that there exists $n\ge k$ so
that $X_{n+1}$ is obtained from $X_n$ by attaching an apartment
$A'$ along $L$. Therefore, the wall $L\subset X$ will be the
intersection of three (actually, four) half-apartments in $X$ and,
hence, thick. \qed

\medskip
This concludes the proof of Theorem \ref{main}.

\begin{rem}
One can modify the construction of the building $X$ as follows.
Besides {\em enlargements}, one allows isometric embeddings $X_n
\to X_{n+1}$, where $X_{n+1}$ is obtained from $X_n$ by attaching
an apartment (or another $X_n'\in \c$) along a convex subcomplex
{\em isometric} to a convex subcomplex in $A$ (resp. in $X_n'$).
\end{rem}

As a corollary of Theorem \ref{main}, we prove

\begin{cor}
Let $\ol{W_{af}}=\R^2 \rtimes W$, the closure of $W_{af}$ in the
full group of isometries of $A$. Then there exists a thick
Euclidean building $X'$ of rank 2 modeled on $(A,\ol{W_{af}})$.
\end{cor}
\proof Let $X$ be a building as in Theorem \ref{main}. Pick a
nonprincipal ultrafilter $\om$ and consider the ultralimit $X_\om$
of the constant sequence $(X)$ with fixed base-point (see e.g.
\cite{KaL} for the definitions). Then, according to \cite{KL},
$X'=X_\om$ is the required building: It is a thick Euclidean
building modeled on $(A,\ol{W_{af}})$. \qed

\section{Concluding remarks}\label{conclusion}

One would like to generalize Theorem \ref{main} in several directions:

1. Eliminate the assumption that $(A, W_{af})$ is special.

2. Allow for uncountable Coxeter groups (besides $\ol{W_{af}}$).

3. Construct buildings $X$ where the automorphism group acts
transitively on the set of apartments.

\medskip
At the moment, we do not know how to deal with 1. However, 2 and 3
can be approached along the same lines. Namely, instead of
constructing a building $X$ as a direct limit of a totally-ordered
(by inclusion) set of spaces $X_n\in \c$, one would like to allow
direct limits using partially ordered sets. This hinges on the
following problem, to which we expect a positive solution:

\begin{prob}
Fix a Coxeter complex $(A, W_{af})$ on which all spaces $X\in \c$
are modeled. Recall that in Section \ref{category}, we defined a
category $\e$ of finite complexes modeled on $(A, W_{af})$, where
morphisms are given by compositions of enlargements.

Find a category $\d$, so that $\e \subset \d \subset \c$, and that
$\d$ admits {\em amalgamations}: Every  diagram
\begin{diagram}
 &  &     X_1  &  & \\
 &\ruTo &  &  & \\
X_0 &  &  &  & \\
 &\rdTo &  &  & \\
 &  &     X_2  &  &
\end{diagram}
extends to a commutative diagram
\begin{diagram}
 &  &     X_1  &  & \\
 & \ruTo &  & \rdTo & \\
X_0 &  & \rTo  &  & X_3\\
 & \rdTo &  & \ruTo & \\
 &  &     X_2  &  &\\
\end{diagram}
It is conceivable that $\d=\e$, i.e., that $\e$ is closed under
amalgamations.
\end{prob}

Defining such $\d$ will solve both 2 and 3: see \cite{Tent} for
the explanation of how to use the amalgamation property to form
buildings with ``large'' automorphism groups. Constructing $\d$ in
the context of spherical buildings modeled on $(S^1, W)$ is the
main technical result of the paper by K.~Tent  \cite{Tent} and is
nontrivial already in this case. Doing this for affine Coxeter
groups appears to be much harder. Existence of the category $\d$
would imply that the affine buildings constructed in our paper
embed in rank 2 affine buildings with highly transitive
automorphism group. A partial confirmation to this comes from the
following observation. Let $X$ be a building constructed as in
Section \ref{mainproof}; such a building is unlikely to have any
nontrivial automorphisms. Let $X'$ be a rank 2 building of
countably infinite thickness which is modeled on $(A, W_{af})$
(i.e., every wall in $X'$ is the intersection of countably
infinitely many distinct half-apartments). Then there exists an
embedding $X\to X'$. In particular, if $W$ is of crystallographic
type, there exists a building $X'$ corresponding to a reductive
algebraic group and, hence, having highly transitive automorphism
group.

Addresses:

\noindent Arkady Berenstein: Department of Mathematics, University
of Oregon, Eugene, OR 97403, USA. (arkadiy@uoregon.edu)

\medskip
\noindent Michael Kapovich: Department of Mathematics, University
of California, Davis, CA 95616, USA. (kapovich@math.ucdavis.edu)


\begin{thebibliography}{BaBE}
\addcontentsline{toc}{section}{Bibliography}


\bibitem[BH]{BH}
M.\ Bridson and A.\ Haefliger,
{\em Metric spaces of non-positive curvature},
Grundlehren, vol. 319, Springer, 1999.






\bibitem[CL]{CL}
R. Charney and A. Lytchak, {\em Metric characterizations of
spherical and Euclidean buildings}, Geom. Topol., vol. 5  (2001),
p. 521--550.

\bibitem[FS]{FS}
M. Funk and K. Strambach, {\em Free constructions}, In:  ``Handbook of incidence geometry,''
739--780, North-Holland, Amsterdam, 1995.


\bibitem[G]{G}
W. Goldman, {\em Geometric structures on manifolds and varieties
of representations}, In  ``Geometry of group representations
(Boulder, CO, 1987),''  p. 169--198, Contemp. Math., vol. 74,
Amer. Math. Soc., 1988.


\bibitem[HKW]{HKW}
P. Hitzelberger, L. Kramer and R. Weiss, {\em Non-discrete
Euclidean buildings for the Ree and Suzuki groups,} Preprint,
2008.

\bibitem[KaL]{KaL}
M. Kapovich and B. Leeb, {\em Asymptotic cones and quasi-isometry
classes of fundamental groups of 3-manifolds}, Geometric Analysis
and Functional Analysis, vol. 5, no. 3 (1995), p. 582--603.

\bibitem[KL]{KL}
B. Kleiner and B. Leeb,
{\em Rigidity of quasi-isometries for symmetric spaces and
 {E}uclidean buildings}, Publ. Math. Inst. Hautes \'Etudes Sci., vol. 86 (1997)
 p. 115--197.

\bibitem[N]{N}
I. Nikolaev, {\em The tangent cone of an Aleksandrov space of
curvature $\leq K$}, Manuscripta Math., vol. 86 (1995), no. 2, p.
137--147.


\bibitem[P]{Parreau}
A.\ Parreau, {\em Immeubles affines: construction par les normes
et \'etude des isom\'etries}, In: ``Crystallographic groups and
their generalizations'' (Kortrijk, 1999), p. 263--302, Contemp.
Math., vol. 262, Amer. Math. Soc., Providence, RI, 2000.

\bibitem[R1]{Ronan1}
M. Ronan, {\em A construction of buildings with no rank 3 residues
of spherical type}, In: ``Buildings and the geometry of diagrams,
(Como, 1984)'', p. 242--248.  Lecture Notes in Math., vol. 1181,
Springer, Berlin, 1986.



\bibitem[R2]{Ronan2}
M. Ronan, ``Lectures on buildings'',  Perspectives in Mathematics, vol. 7,
Academic Press, 1989.

\bibitem[Te]{Tent}
K. Tent,  {\em Very homogeneous generalized $n$-gons of finite
Morley rank}, J. London Math. Soc. (2), vol. 62 (2000), no. 1, p.
1--15.

\bibitem[T1]{Tits1}
J. Tits, ``Buildings of spherical types and finite BN-pairs'',
Lecture Notes in Math., vol. 386, Springer, Berlin, 1974.

\bibitem[T2]{Tits2}
J. Tits, {\em Endliche Spiegelungsgruppen, die als Weylgruppen
auftreten}, Invent. Math., vol. 43 (1977), no. 3, p. 283--295.

\bibitem[T3]{Tits3}
J. Tits, {\em Immeubles de type affine},  In ``Buildings and the
geometry of diagrams (Como, 1984)'', p. 159--190.  Lecture Notes
in Math., vol. 1181, Springer, Berlin, 1986.


\end{thebibliography}
\end{document}